\documentclass[psfig, 12pt]{article}
\usepackage{hyperref,amssymb,amsmath,graphicx,verbatim,amsfonts,amsthm}
\usepackage{psfig}
\newtheorem{thm}{Theorem}[section]
\newtheorem{cor}[thm]{Corollary}
\newtheorem{lem}[thm]{Lemma}
\newtheorem{prop}[thm]{Proposition}
\newtheorem{clm}[thm]{Claim}

\newtheorem{conj}[thm]{Conjecture}

\newtheorem*{thm*}{Theorem}

\theoremstyle{definition}
\newtheorem{dfn}[thm]{Definition}
\theoremstyle{remark}

\numberwithin{equation}{section}

\newcommand{\Claim}{\textbf{Claim.}\hspace{10pt}}
\newcommand{\Rem}{\textbf{Remark.}\hspace{10pt}}

\newcommand{\norm}[1]{\left\Vert#1\right\Vert}
\newcommand{\abs}[1]{\left\vert#1\right\vert}

\newcommand{\IP}[1]{\left<#1\right>}

\newcommand{\set}[1]{\left\{#1\right\}}
\newcommand{\Set}[2]{ \left\{#1 \ \big| \ #2 \right\} }
\newcommand{\br}[1]{\left[#1\right]}
\newcommand{\Br}[2]{ \left[#1 \ \big| \ #2 \right] }
\newcommand{\sr}[1]{\left(#1\right)}

\newcommand{\Integer}{\mathbb{Z}}

\newcommand{\Z}{\Integer}
\newcommand{\N}{\mathbb{N}}

\newcommand{\R}{\mathbb{R}}

\newcommand{\eps}{\varepsilon}

\newcommand{\E}{\mathbf{E}}
\newcommand{\Var}{\mathbf{Var}}

\newcommand{\1}[1]{\mathbf{1}_{\set{ #1 } }}
\newcommand{\ov}[1]{\overline{#1}}

\newcommand{\eqdef}{\stackrel{\textrm{\tiny def} }{=}}

\def\squareforqed{\hbox{\rlap{$\sqcap$}$\sqcup$}}
\def\qed{\ifmmode\squareforqed\else{\unskip\nobreak\hfil
\penalty50\hskip1em\null\nobreak\hfil\squareforqed
\parfillskip=0pt\finalhyphendemerits=0\endgraf}\fi}

\newcommand{\GCD}{$G$-Cylinder-DLA }

\newcommand{\z}{\zeta}
\renewcommand{\k}{\kappa}

\newcommand{\eqn}[1]{ \begin{eqnarray*}
  #1
\end{eqnarray*} }

\newcommand{\m}{\mathfrak{m}}

\newcommand{\X}{\Xi}
\begin{document}

\title{Diffusion Limited Aggregation on a Cylinder }

\author{Itai Benjamini \thanks{Email:
\texttt{itai.benjamini@weizmann.ac.il}. } \and Ariel Yadin
\thanks{Email: \texttt{ariel.yadin@weizmann.ac.il}. } }

\date{ }

\maketitle

\begin{abstract}
We consider the DLA process  on a cylinder $G \times \N$. It is
shown that this process ``grows arms'', provided that the base
graph $G$ has small enough mixing time. Specifically, if the
mixing time of $G$ is at most $\log^{(2-\eps)}\abs{G}$, the time
it takes the cluster to reach the $m$-th layer of the cylinder is
at most of order $m \cdot \frac{ \abs{G}}{\log\log \abs{G}}$. In
particular we get examples of infinite Cayley graphs of degree
$5$, for which the DLA cluster on these graphs has arbitrarily
small density.

In addition, we provide an upper bound on the rate at which the
``arms'' grow. This bound is valid for a large class of base
graphs $G$, including  discrete tori of dimension at least $3$.

It is  also shown that for any base graph $G$, the density of the
DLA process on a $G$-cylinder is related to the rate at which the
arms of the cluster grow.  This implies, that for any vertex
transitive $G$, the density of DLA on a $G$-cylinder is bounded by
$2/3$.
\end{abstract}

\section{Introduction}

Diffusion Limited Aggregation (DLA), is a growth model introduced
by Witten and Sander (\cite{WittenSander}). The process starts
with a particle at the origin of $\Z^d$.  At each time step, a new
particle starts a simple random walk on $\Z^d$ from infinity (far
away). The particle is conditioned to hit the existing cluster
(when $d \geq 3$).  When the particle first hits the outer
boundary of the cluster, it sticks and the next step starts,
forming a growing family of clusters.

We consider a variant of this model, where the underlying graph of
the process is a cylinder with base $G$, $G$ being some finite
graph. A precise definition is given in Section \ref{scn:
Definition of CDLA process}.

This paper contains three main results:

The first, Theorem \ref{thm:  E [ T_m ] leq o(nm)}, states that if
$G$ has small enough mixing time, then the time it takes the
cluster to reach the $m$-th layer of the cylinder is $o(m \cdot
\abs{G})$, were $\abs{G}$ is the size of $G$.  In fact, for a
graph $G$ with mixing time at most $\log^{(2-\eps)} \abs{G}$ (for
any constant $\eps$), the time to reach the $m$-th layer is at
most of order $m \cdot \frac{ \abs{G}}{ \log\log \abs{G} }$. This
phenomenon is sometimes dubbed as ``the aggregate grows arms'',
i.e. grows faster than order $\abs{G}$ particles per layer. The
analogous phenomenon in the original DLA model on $\Z^d$ is
considered a notoriously difficult open problem.  In
\cite{Kesten1, Kesten2, Kesten3}, Kesten provides upper bounds on
the growth rate of the DLA aggregate in $\Z^d$. Eberz-Wagner
\cite{Eberz} proved the existence of infinitely many holes in the
two-dimensional DLA aggregate.

The second result concerns the density of the limit cluster, the
union of all clusters obtained at some finite time. Theorem
\ref{thm: E[D] leq lim (1/mn) E[T_m']} shows that the expected
rate at which the cluster grows bounds this density. This has two
implications:
\begin{enumerate}

    \item Theorem \ref{thm:  for any v.t. graph D leq 2/3} states
    that for any vertex transitive graph $G$, the DLA process on the
    $G$-cylinder has density bounded by $2/3$. This includes the
    cases where $G$ is a $d$-dimensional Torus.

    \item Theorem \ref{thm: D(n) leq 3/log log n for RM family of
    graphs} shows that for $G$ with small enough mixing time, the
    density tends to $0$ as the size of $G$ tends to infinity.

\end{enumerate}

Finally, Theorem \ref{thm:  lower bound on speed} is a lower bound
on the expected time the cluster reaches the $m$-th layer,
complementing the upper bound in Theorem \ref{thm: E [ T_m ] leq
o(nm)}.  This lower bound implies that the cluster cannot grow too
fast, and in fact for many natural graphs it cannot grow faster
than $\abs{G}^{c}$ for some universal $0 < c < 1$. The lower bound
holds for a wider range of graphs at the base of the cylinder than
the upper bound (including $d$-dimensional tori for $d \geq 3$).

We remark that our estimates for the upper bound are crude, and
simulations indicate that there is much room for improvement. In
fact we believe the truth to be closer to the lower bound, see
Conjecture \ref{conj: Conjecture about T_m}. Proving Conjecture
\ref{conj: Conjecture about T_m} will imply that for \emph{any}
family of graphs $\set{G_n}$, the density of the DLA process on
the $G_n$-cylinder tends to $0$ as the size of $G_n$ tends to
infinity (see Conjecture \ref{conj: Conjecture about T_m}).

For other very different variants of one dimensional DLA see
\cite{1dimDLA, KestenVladas}.  Another paper dealing with
random-walk related questions on cylinders with varying bases is
\cite{Alain}.

The rest of this paper is organized as follows:

First we introduce some notation.  In Section \ref{scn: Definition
of CDLA process} we define the process, and random variables
associated with it.  In Section \ref{scn: Statement of Results} we
state the first main result.  Section \ref{scn: Main tool for
proving theorems} is devoted to proving Theorem \ref{thm:  E [
T(M+1) - t ] leq  n/loglog n}, the main tool used to prove the
main results. After the formulation of this theorem, a sketch of
the key dichotomy idea is given, followed by a short discussion.
In Section \ref{scn: Density}, we define the density of the DLA
process on a cylinder. We also prove the theorems bounding the
density in the above mentioned cases, Theorems \ref{thm:  for any
v.t. graph D leq 2/3} and \ref{thm: D(n) leq 3/log log n for RM
family of graphs}. Finally, in Section \ref{scn:  lower bound
section} we prove the lower bound on the growth rate of the
cluster, Theorem \ref{thm: lower bound on speed}.

Let us note that the set up of DLA on a cylinder suggests another
natural problem we are now pursuing. That is, how long does it
take until the cluster clogs the cylinder? (This problem may be
related to \cite{DemboSznit}.)

Other possible directions for further research are presented in
the last section, followed by an appendix which contains a few
standard variants on some simple random walk results we need.

\paragraph{Acknowledgement.} We wish to thank Amir Yehudayoff for
many useful discussions, and remarks about a preliminary version
of this note.

\subsection{Notation}

Let $G$ be a graph.  $V(G)$ and $E(G)$ denote the vertex set and
edge set of $G$ respectively. We use the notation $v \in G$ to
denote $v \in V(G)$. For two vertices $u,v$ in $G$ we use the
notation $u \sim v$ to denote that $u$ and $v$ are adjacent.

%

For a graph $G$, define the \emph{cylinder with base $G$}, denoted
$G \times \N$, by:  The vertex set of $G \times \N$, is the set
$V(G) \times \N$.  The edge set is defined by the following
relations: For all $u,v \in G$ and $m,k \in \N$, $(u,m) \sim
(v,k)$ if and only if: either $m=k$ and $u \sim v$, or $\abs{m-k}
= 1$ and $u=v$.  The cylinder with base $G$ is just placing
infinitely many copies of $G$ one over the other, and connecting
each vertex in a copy to its corresponding vertices in the
adjacent copies.

By the simple random walk on a graph, we refer to the process
where at each step the particle chooses a neighbor uniformly at
random and moves to that neighbor.  By the lazy random walk with
holding probability $\alpha$, we mean a walk that with probability
$\alpha$ stays at its current vertex, and with probability
$1-\alpha$ chooses a neighbor uniformly at random.  By lazy random
walk (without stating the holding probability) we refer to the
walk that chooses uniformly at random from the set of neighbors
and the current vertex. A lazy random walk is a simple random walk
on the same graph with a self loop added at each vertex.

For simplicity, this paper will only deal with regular graphs;
i.e. graphs such that all vertices are of the same degree.

We define the notion of the \emph{mixing time} of a $d$-regular
graph $G$: Let $\set{g'_t}_{t \geq 0}$ be a lazy random walk on
$G$. The mixing time of $G$, is defined by
\begin{equation}
\m = \m(G) \eqdef \min \Set{ t > 0 }{ \forall \ u,v \in G \
\forall \ s \geq t \ , \ \Pr \Br{ g'_s = u }{ g'_0 = v } \geq
\frac{1}{2 \abs{G}} } .
\end{equation}
This is a valid definition, since for all $v \in G$,
$$ \lim_{t \to \infty} \max_{u \in G} \abs{ \Pr \Br{ g'_t = u }{ g'_0 = v } -
\frac{1}{\abs{G}} } = 0 . $$ %
(This can be seen via Lemma \ref{lem:  mixing time of RW}.)

For a probability event $A$, we denote by $\ov{A}$ the complement
of $A$.





\section{Cylinder DLA} \label{scn:  Definition of CDLA process}

\subsection{Definition}

Fix a graph $G$. We define the \emph{$G$-Cylinder-DLA} process:

Consider the graph $G \times \N$.  Denote by $G_i$ the induced
subgraph on the vertices $V(G) \times \set{i}$, for all $i \in
\N$. We call $G_i$ the \emph{$i$-th layer} of $G \times \N$.

The process is an increasing sequence, $\set{A_t}_{t=0}^\infty$,
of connected subsets of $G \times \N$.  We start with $A_0 = G_0$.
Given $A_t$, define the set $A_{t+1}$ as follows:

Let $\partial A_t$ be the set of all vertices of $G \times \N$
that are not in $A_t$, but are adjacent to some vertex of $A_t$.
That is,
$$ \partial A_t = \Set{ u \in G \times \N }{ u \not\in A_t, \
\exists v \in A_t \ : \ u \sim v } . $$

Let a particle perform a simple random walk on $G \times \N$
starting from infinity, and stop when the particle hits $\partial
A_t$.  Let $u$ be the vertex in $\partial A_t$ where the particle
is stopped.  Then, set $A_{t+1} = A_t \cup \set{u}$.

We find it convenient to use the following alternative (but
equivalent) definition:

Let $M(t)= \min \Set{ i \in \N}{ G_i \cap A_t = \emptyset }$. That
is, $M(t)$ is the lowest layer of $G \times \N$ that does not
intersect the cluster $A_t$. Let $(g_{t+1}(i), \z_{t+1}(i)) \in G
\times \N$, $i=0,1,2,\ldots$, be a simple random walk on $G \times
\N$, such that $g_{t+1}(0)$ is uniformly distributed in $G$, and
$\z_{t+1}(0) = M(t)$.

Let $\k(t+1)$ be the first time at which the walk is in $\partial
A_t$. That is,
$$ \k(t+1) = \min \Set{ r \geq 0 }{ (g_{t+1}(r), \z_{t+1}(r)) \in \partial A_t } . $$
Since the walk is recurrent, $\k(t+1) < \infty$ with probability
$1$. Let $\k= \k(t+1)$.  Then, $(g_{t+1}(\k), \z_{t+1}(\k) )$ is
distributed on the set $\partial A_t$. Set $A_{t+1} = A_t \cup
\set{ (g_{t+1}(\k), \z_{t+1}(\k) ) }$.

This construction is equivalent to ``starting from infinity''; a
simple random walk starting at higher and higher layers, will take
more and more steps before reaching the layer $M(t)$.  Thus, as
the starting layer tends to infinity, the distribution of the
particle at the first time it hits the layer $M(t)$ is tending to
uniform.



\begin{figure}[htbp]

\centerline{ \psfig{figure=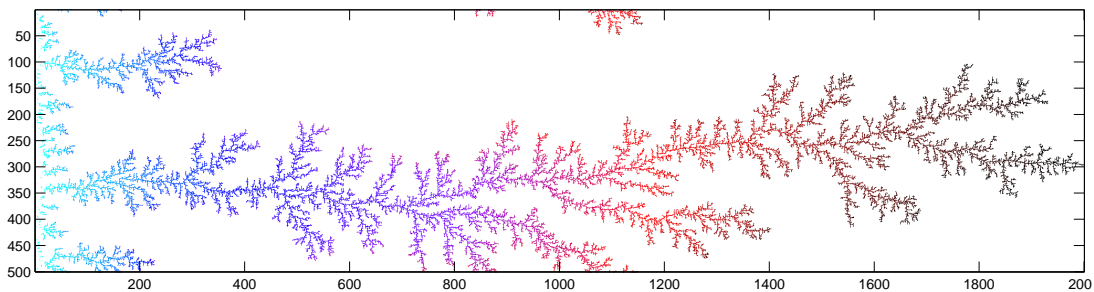,width=15cm,clip=} }

\caption{$G$-Cylinder-DLA, where $G$ is the cycle on $500$
vertices. The number of particles is approximately $64,400$.}
\label{figure: CDLA-cycle}
\end{figure}


\subsection{ }

Let $\set{A_t}_{t=0}^{\infty}$ be a \GCD process.  $A_t$ is called
the (\GCD) \emph{cluster at time $t$}.  Define the following
random variables:

For $A_t$, the \GCD cluster at time $t$, define the \emph{load of
the $i$-th layer} at time $t$ by
$$ L_t(i) = \abs{ A_t \cap G_i } . $$
$L_t(i)$ is the number of particles in the cluster at time $t$, on
the $i$-th layer.  Also define
$$ L_t(\geq i) = \sum_{j \geq i} L_t(j) , \quad \textrm{ and } \quad
L_t(> i) = \sum_{j > i} L_t(j) . $$ %
$L_t(\geq i)$ (respectively $L_t(> i)$) is the total load on
layers $\geq i$ (respectively $>i$). Note that
$$ L_t(\geq i) = \sum_{j = i}^{M(t)-1} L_t(j) , \quad \textrm{ and } \quad
L_t(> i) = \sum_{j = i+1}^{M(t)-1} L_t(j) . $$ %
When subscripts become too small, we write $L(t,i)$ instead of
$L_t(i)$ (and similarly for ${L(t,\geq i)}$ and ${L(t,>i)}$).


Here are some properties of the Cylinder-DLA process, that we
leave for the reader to verify.  (This can help to get used to the
notation.)
\begin{enumerate}

    \item For all $s > t$, $A_t \subsetneq A_s$.
    \item For any $i \in \N$, and $s \geq t$, $L_t(i) \leq L_s(i)$.
    \item If $L_t(i) = 0$, then $L_t(j) = 0$ for all $j \geq i$.
    \item For all $i \geq M(t)$, $L_t(i) = 0$.  For all $i <
    M(t)$, $L_t(i) \geq 1$.
    \item For all $t$,
    $$ \sum_{i=0}^{\infty} L_t(i) = \sum_{i=0}^{M(t)-1} L_t(i) .
    $$ %
    \item The following events are identical (for any $t > 0$):
    $$ \set{ L_t(\geq i) > L_{t-1}(\geq i) } = \set{ L_t(\geq i) =
    1 + L_{t-1}(\geq i) } = \set{ \z_t(\k(t)) \geq i } . $$
\end{enumerate}

\subsection{\GCD grows arms, for quickly mixing $G$} \label{scn: Statement of
Results}

\begin{thm} \label{thm:  E [ T_m ] leq o(nm)}
Let $2 \leq d \in \N$. There exists $n_0 = n_0(d)$, such that the
following holds for all $n > n_0$:

Let $G$ be a $d$-regular graph of size $n$, and mixing time
$$ \m(G) \leq \frac{\log^2(n)}{(\log \log (n))^5} . $$ %
Let $\set{A_t}$ be a \GCD process. For $m \in \N$, define
$$ T_m = \min \Set{ t \geq 0 }{ A_t \cap G_m \neq \emptyset } . $$
$T_m$ is the time the cluster first reaches the layer $m$.

Then, for all $m$,
$$ \E \br{ T_m} < \frac{4 m n}{\log \log n} . $$
\end{thm}

The proof of Theorem \ref{thm:  E [ T_m ] leq o(nm)} is via
Theorem \ref{thm:  E [ T(M+1) - t ] leq n/loglog n} below.

\Rem One may suggest that the reason Theorem \ref{thm:  E [ T_m ]
leq o(nm)} can be proved, is that we use for the base graph $G$,
graphs that are so highly connected that in some sense there is no
geometry. We stress that the class of graphs that have
$\log^{(2-\eps)} \abs{G}$ mixing time, is much larger than what is
known as ``expander graphs''.  This class includes many natural
families of graphs, including lamplighter graphs on tori of
dimension $2$ and above (see \cite{PeresRevelle}).

We remark that Theorem \ref{thm:  E [ T(M+1) - t ] leq  n/loglog
n} below is in some sense a ``worst case'' analysis.  Thus, we
believe that our results are not optimal.  In fact, we conjecture
that a stronger result than Theorem \ref{thm:  E [ T_m ] leq
o(nm)} should hold for any graph at the base of the cylinder:

\begin{conj} \label{conj: Conjecture about T_m}
Let $\set{G_n}$ be a family of $d$-regular graphs such that
$\lim_{n \to \infty} \abs{G_n} = \infty$. There exist $0< \gamma <
1$ and $n_0$ such that for all $n > n_0$ the following holds:

Set $G = G_n$ and let $\set{A_t}$ be a \GCD process. For $m \in
\N$, define
$$ T_m = \min \Set{ t \geq 0 }{ A_t \cap G_m \neq \emptyset } . $$
$T_m$ is the time the cluster first reaches the layer $m$.

Then, for all $m$,
$$ \E \br{ T_m} \leq m \abs{G_n}^{\gamma} . $$
\end{conj}


\section{The time to stick to a new layer} \label{scn: Main tool
for proving theorems}

The following theorem states that under the assumption that $G$
has small enough mixing time, in the $G$-Cylinder-DLA process, the
expected amount of particles until one sticks to the new layer, is
substantially less than $\abs{G}$.

Note that since for all $m$, we can write the telescopic sum
$$ T_m = \sum_{\ell=1}^m (T_\ell - T_{\ell-1}) , $$
Theorem \ref{thm:  E [ T_m ] leq o(nm)} follows from Theorem
\ref{thm:  E [ T(M+1) - t ] leq  n/loglog n}, by linearity of
expectation.

\begin{thm} \label{thm:  E [ T(M+1) - t ] leq  n/loglog n}
Let $2 \leq d \in \N$. There exists $n_0 = n_0(d)$, such that the
following holds for all $n > n_0$:

Let $G$ be a $d$-regular graph of size $n$, and mixing time
$$ \m(G) \leq \frac{\log^2(n)}{(\log \log (n))^5} . $$
Let $A_t$ be a \GCD cluster at time $t$. Define
$$ T = \min \Set{ s > t }{ M(s) > M(t) } . $$
$T$ is the first time that a particle sticks to the empty layer,
$G_{M(t)}$. Then,
$$ \E \br{ T - t } < \frac{4 n}{\log \log n} . $$
\end{thm}

In order to prove Theorem \ref{thm:  E [ T(M+1) - t ] leq n/loglog
n}, we need a few lemmas, stated and proved in this section. The
proof of Theorem \ref{thm:  E [ T(M+1) - t ] leq n/loglog n} is
deferred to Section \ref{scn:  Proof of Theorem E [ T(M+1) - t ]
leq n/loglog n}.  The main idea of the proof is in the following
proof sketch:

\textit{Proof Sketch.  }  The cluster $A_t$ can be in two states:
Either it is such that particles stick quickly to it; i.e.
particles take few steps before sticking to the cluster.  Or, the
particles take many steps before sticking to the cluster.

In the first case, the particles take few steps before sticking.
Thus, the particles cannot stick many layers below $M(t)$, so they
build up a heavy load on the layers near $M(t)$. Each time a layer
has a heavy load, there is better chance of the next particles to
stick to the layers above it.  So, in less than $O \sr{
\frac{n}{\log\log n} }$ particles, there is a heavy load on the
layer $M(t)-1$, and the probability of sticking to the layer
$M(t)$ is now substantially greater than $1/n$. This case is dealt
with in Lemma \ref{lem:  if k(t+1) is small then E[ T] leq
n/loglog(n)}.

%

In the second case, the particles take many steps before sticking.
Thus, they also make many long excursions above the layer $M(t)$.
Because the base of the cylinder, $G$, has small enough mixing
time, after each such excursion, there is a chance of at least
$1/2n$ to stick to the layer $M(t)$.  This occurs many times, so
the probability of sticking to the layer $M(t)$ is much greater
than $1/n$. This case is dealt with in Corollary \ref{cor: when
there are many steps, for RM graphs}.

The proof of Theorem \ref{thm:  E [ T(M+1) - t ] leq n/loglog n}
in Section \ref{scn:  Proof of Theorem E [ T(M+1) - t ] leq
n/loglog n} combines both cases, to show that in both cases, the
expected time until a particle sticks to the new layer $M(t)$, is
substantially smaller than $n$.

\Rem As stated above, the proof of Theorem \ref{thm:  E [ T(M+1) -
t ] leq n/loglog n} is in some sense a ``worst-case'' analysis.
The first part, regarding the case where particles take few steps
before sticking, is valid for any regular $G$ (not only those with
small mixing time). But in reality, simulations show that this is
not what really happens. The particles do not build a series of
higher and higher layers with large loads.

On the other hand, the second part, (where particles take many
steps and thus return to the layer $M(t)$ many times, thus
increasing the probability of sticking to $M(t)$) is probably what
does actually occur.  In fact, we suspect that this is true not
only for graphs with small mixing time, but for any graph at the
base of the cylinder (see Conjecture \ref{conj: Conjecture about
T_m}).

\Rem It may be of use to note that Theorem \ref{thm:  E [ T(M+1) -
t ] leq  n/loglog n} holds also if $A_t$ is replaced with any
subset of $G \times \N$ intersecting all layers up to $M(t)$. In
particular, given \emph{any} cluster, not necessarily grown by a
\GCD process, the expected time until a particle sticks to the new
layer is bounded by order $\frac{\abs{G} }{\log\log \abs{G}}$.

\subsection{A large load on a high layer}

In this section, we show that if there is a high enough layer
($\geq M(t) - \frac{\log n}{4 \log \log n}$) with large load (at
least $\frac{n}{\log n}$), then the expected time until a particle
sticks to the new layer, $M(t)$, is $o(n)$.


\begin{lem} \label{lem:  Expected time to reach new layer, if
there is a wall} %
There exists $n_0$, such that the following holds for all $n >
n_0$: 
Let $G$ be a $d$-regular graph of size $n$. Set
$$ \mu = \mu(n) = \left\lfloor \frac{ \log (n)}{4 \log \log (n)} \right\rfloor ,
\quad \textrm{ and } \quad  \nu = \nu(n) = \log(n) . $$ %
Let $A_t$ be a \GCD cluster at time $t$.  Let
$$ T = \min \Set{ s > t }{ M(s) > M(t) } . $$
$T$ is the first time that the cluster reaches the new layer.

Assume that there exists $j \geq M(t)-\mu$ such that $L_t(j) \geq
\frac{n}{\nu}$.  Then,
$$ \E \br{ T-t } \leq \frac{n}{ 4 \log \log (n) } . $$
\end{lem}

The main idea of the proof is as follows:  If a layer $j$ has load
$m$, then the probability to stick above layer $j$ is at least
$m/n$.  Thus, to get a layer $i > j$ with load
$\frac{m}{\log(m)}$, we need $o(n)$ particles.  Thus, building
higher and higher layers with high loads, we reach the empty layer
in $o(n)$ particles.

\begin{proof}[Proof of Lemma \ref{lem:  Expected time to reach new layer, if
there is a wall}] %
The following proposition states that if there is a layer with
load $m$, then the probability of particles sticking above that
layer is at least $m/n$.

\begin{prop} \label{prop:  Pr [ I_s ] geq m/n}
Let $G$ be a $d$-regular graph of size $n$. Let $A_t$ be a \GCD
cluster at time $t$. Fix a layer $j>0$. Assume that $L_t(j) \geq
m$.  For $s > t$, let $I_s$ be the indicator function of the event
that the $s$-th particle sticks to a layer $\geq j+1$.  That is,
$$ I_s = \1{ L(s,\geq j+1) > L(s-1,\geq j+1) } . $$
Then, for all $s > t$,
$$ \Pr \Br{ I_s = 1 }{ I_r \ , \ t < r<s  } \geq \frac{m}{n} , $$
for any values of $I_r$, $t<r<s$.
\end{prop}

\begin{proof}
Set $s>t$.  Condition on the values of $I_r$, $t<r<s$.  Let
$A_{s-1}$ be the cluster at time $s-1$.  Let $(g(\cdot),\z(\cdot))
= ((g_s(\cdot),\z_s(\cdot))$ be the walk of the $s$-th particle.
Note that for any $r$, if $(g(r),\z(r)) \in \partial A_{s-1} \cup
A_{s-1}$, then $\k(s) \leq r$.  At time $s-1$ the layer $j$ has
load $L_{s-1}(j) \geq L_t(j) \geq m$.  Thus,
$$ \abs{ \sr{ \partial A_{s-1} \cup A_{s-1} } \cap G_{j+1} } \geq
m . $$ %
Let $k$ be the first time the walk $((g(\cdot),\z(\cdot))$ hits
the layer $j+1$. Then, since the uniform distribution on $G$ is
the stationary distribution, $g(k)$ is uniformly distributed in
$G_{j+1}$.  Thus, %
\eqn{ \Pr \Br{ \k(s) \leq k }{ A_{s-1} } & \geq & \Pr \Br{
(g(k),\z(k)) \in \partial A_{s-1} \cup A_{s-1} }{ A_{s-1} } \geq \frac{m}{n} . } %
Since for all $0 \leq r \leq k$ we have that $\z(r) \geq j+1$, we
get that
$$ \Pr \Br{ I_s = 1 }{ A_{s-1} } \geq \Pr \Br{ \k(s) \leq k }{
A_{s-1} } \geq \frac{m}{n} . $$ %
Let ${\cal A}$ be the set of all clusters $A \subseteq G \times
\N$ such that $\Pr \Br{ A_{s-1} = A }{ I_r \ , \ t<r<s } > 0$.
Then we have,
\begin{eqnarray*}
    & & \Pr \Br{ I_s = 1 }{ I_r \ , \ t < r<s  } \\
    & = & \sum_{A \in {\cal A}} \Pr \Br{ I_s = 1 }{ A_{s-1} = A }
    \Pr \Br{ A_{s-1} = A }{ I_r \ , \ t<r<s } \geq \frac{m}{n} .
\end{eqnarray*}
\end{proof}

Assume there is a layer with load $m$.  Since each particle sticks
above this layer with probability at least $m/n$, the expected
time until there are $\ell$ new particles above this layer should
be at most $\ell \cdot (n/m)$.  This is captured in the following
proposition:

\begin{prop} \label{prop: Expected time for ell particles above wall}
Let $G$ be a $d$-regular graph of size $n$. Let $A_t$ be a \GCD
cluster at time $t$. Fix a layer $j > 0$. Assume that $L_t(j) \geq
m$. For $\ell \in \N$, define
\begin{equation} \label{eqn:  Dfn of S(ell)}
S_{\ell} = \min \Set{ s \geq t }{ L_s(\geq j+1) = \ell + L_t(\geq
j+1) } .
\end{equation}
That is, $S_{\ell}$ is the first time that there are $\ell$ new
particles in the layers $\geq j+1$ (so $S_0 = t$).

Then,
$$ \E \br{ S_{\ell} - t} = \E \br{ S_{\ell} - S_0} \leq \ell \frac{n}{m} . $$
\end{prop}

\begin{proof}
By Proposition \ref{prop:  Pr [ I_s ] geq m/n}, for all $k \geq
1$, $S_{k} - S_{k-1}$ is dominated by a geometric random variable
with mean $\leq \frac{n}{m}$.  Thus,
$$ \E \br{ S_{\ell} - t} = \sum_{k=1}^{\ell} \E \br{ S_{k} -
S_{k-1} } \leq \ell \frac{m}{n} . $$
\end{proof}

%
%

With these two propositions,  we continue with the proof of Lemma
\ref{lem: Expected time to reach new layer, if there is a wall}.

Set $M = M(t)$. Set $T_0 = t$. For $r \geq 0$, define inductively
the following stopping times:
$$ T_r = \min \Set{ s \geq T_{r-1} }{ \exists i \geq j+r \ : \ L_s(i)
\geq n \nu^{-(2r+1)} } . $$ %
That is, $T_r$ is the first time that there exists a ``high
enough'' layer (higher than $j+r$), such that the load on that
layer is ``large enough'' (larger than $n \nu^{-(2r+1)}$).

Consider time $T_{\mu}$.  At this time, we have that there exists
a layer $i \geq j+\mu \geq M$ such that $L_{T_\mu}(i) \geq n
\nu^{-(2\mu+1)} \geq 1$. So $M(T_{\mu}) > M$ and $T \leq T_\mu$.
Thus, we can write
$$ T - t = \sum_{r=1}^\mu \sr{ \min \set{ T, T_r } - \min \set{
T, T_{r-1} } } . $$

For all $r \geq 0$, set $\tau(r) = \min \set{ T, T_r }$.

\Claim For all $r > 0$,
$$ \E \br{ \tau(r) - \tau(r-1) } \leq \frac{n}{\nu} . $$

\begin{proof}
Fix $r>0$.  For $\ell \in \N$, define
$$ S_{\ell} = \min \Set{ s \geq T_{r-1} }{ L_s(\geq j+r) = \ell +
L_{T_{r-1}}(\geq j+r) } . $$ %
That is, $S_{\ell}$ is the first time that there are $\ell$ new
particles in the layers $\geq j+r$. 
So $S_0 = T_{r-1}$. Let $a = \mu \lceil n \nu^{-(2r+1)} \rceil $.

\begin{itemize}
\item[Case 1:] $T \leq T_{r-1}$. Then $\tau(r) - \tau(r-1) = 0 \leq S_a -
S_0$.

\item[Case 2:] $T > T_{r-1}$ and $T \leq S_a$. Then $\tau(r) -
\tau(r-1) \leq T - T_{r-1} \leq S_a - S_0$.

\item[Case 3:] $T > T_{r-1}$ and $T > S_a$. Note that if $T > S_a$,
then $M(S_a) = M$. At time $S_a$, there are at least $a$ particles
on the layers $\geq j+r$. So, if $T > S_a$ then
$$ a \leq \sum_{i=j+r}^{M(S_a)-1} L_{S_a}(i) =
\sum_{i=j+r}^{M-1} L_{S_a}(i) . $$ %
So there exists some $j+r \leq i \leq M-1$ such that $L_{S_a}(i)
\geq \frac{a}{M-(j+r)}$. Since $j \geq M-\mu$, we have $L_{S_a}(i)
\geq \frac{a}{\mu} \geq n \nu^{-(2r+1)}$.

So we conclude that if $T > S_a$ then $T_r \leq S_a < T$.  So $T >
S_a$ implies that
$$ \tau(r) - \tau(r-1) = T_r - T_{r-1} \leq S_a - S_0 . $$
\end{itemize}
Thus, in all three cases, $\tau(r) - \tau(r-1) \leq S_a - S_0$.

At time $S_0 = T_{r-1}$, by the definition of $T_{r-1}$, we have
that for some $i \geq j+r-1$, there is a load $L_{S_0}(i) \geq n
\nu^{-(2r-1)}$ (for $r=1$ we can choose $i=j$, and since $T_0 = t$
we have by assumption that $L_t(j) \geq \frac{n}{\nu}$). By
Proposition \ref{prop: Expected time for ell particles above
wall}, with $j=i, t=T_{r-1} = S_0$ and $m = n \nu^{-(2r-1)}$, we
have that for large enough $n$
$$ \E \br{ \tau(r) - \tau(r-1) } \leq \E \br{ S_a - S_0 } \leq
a \frac{n}{n \nu^{-(2r-1)} } \leq \frac{n}{\nu} . $$
\end{proof}

Returning to the proof of Lemma \ref{lem:  Expected time to reach
new layer, if there is a wall}, for all $r > 0$,
$$ \E \br{ \tau(r) - \tau(r-1) } \leq \frac{n}{\nu} . $$
Thus, for large enough $n$,
$$ \E \br{ T - t } = \sum_{r=1}^\mu \E \br{ \tau(r) - \tau(r-1)}
\leq \mu \cdot \frac{n}{\nu} \leq \frac{n}{ 4 \log\log (n) } . $$
\end{proof}


\subsection{Particles take few steps}

Recall that $\k(s)$ is the number of steps the $s$-th particle
takes until it sticks (to $\partial A_{s-1}$).  In this section,
we show that if $\k(t+1)$ is small, then all particles $s>t$, have
a good chance of sticking at high layers.  Thus, a small amount of
particles is needed to get a high layer with large load.

\begin{lem} \label{lem:  if k(t+1) is small then E[ T] leq
n/loglog(n)} %
There exists $n_0$ such that the following holds for all $n >
n_0$: Let $G$ be a $d$-regular graph of size $n$. Set
$$ \mu = \mu(n) = \left\lfloor \frac{ \log (n)}{4 \log \log (n)} \right\rfloor ,
\quad \textrm{ and } \quad  \nu = \nu(n) = \log(n) . $$ %
Let $A_t$ be a \GCD cluster at time $t$.  Let
$$ T = \min \Set{ s > t }{ M(s) > M(t) } . $$
$T$ is the first time that the cluster reaches the new layer.

Assume that $\Pr \br{ \k(t+1) \leq \frac{\mu^2}{4} } \geq
\frac{1}{4}$. Then,
$$ \E \br{ T-t } \leq \frac{5 n}{2 \log \log (n)} . $$
\end{lem}

\begin{proof}

In the following two propositions, we use the fact that with
probability at least $1/4$, the particle takes a small amount of
steps to stick.

\begin{prop} \label{prop:  the particle sticks without going below the layer M(t)-x} %
Let $G$ be a $d$-regular graph. Let $A_{t}$ be a \GCD cluster at
time $t$, and consider the $(t+1)$-th particle.  Let $y \in \N$
and assume that $\Pr \br{ \k(t+1) \leq y^2/4 } \geq \frac{1}{4}$.
Then,
$$ \Pr \br{ \min_{0 \leq r \leq \k(t+1)} \z_{t+1}(r) \geq M(t)- y } \geq \frac{1}{8}
. $$ %
That is, with probability at least $1/8$, the particle sticks
without ever going below the layer $M(t)- y$.
\end{prop}

\begin{proof}
Set $x = \lfloor \frac{y^2}{4} \rfloor$.  Note that
$$ \Pr \br{ \k(t+1) \leq x } \geq \frac{1}{4} . $$

Let $(g(\cdot),\z(\cdot)) = (g_{t+1}(\cdot),\z_{t+1}(\cdot))$ be
the walk the $(t+1)$-th particle takes. That is, $g(0)$ is
uniformly distributed in $G$, and $\z(0) = M(t)$.  Let $\k =
\k(t+1)$ be the first time the walk hits $\partial A_{t}$.

Note that
$$ \set{ \min_{0 \leq r \leq x} \z(r) \geq M(t) - y } \cap
\set{ \k \leq x } \quad \textrm{ implies } \quad \set{ \min_{0
\leq r \leq \k} \z(r) \geq M(t) - y } . $$

The walk $\z(0), \ldots, \z(x)$, is an $x$-step lazy random walk,
with holding probability $1-\alpha = \frac{d}{d+2}$. By Lemma
\ref{lem:  Lazy RW does not pass x in x^2 steps}, we have that
$$ \Pr \br{ \min_{0 \leq r \leq x} \z(r) \geq M(t) - y } \geq
\Pr \br{ \max_{1 \leq r \leq x } \abs{\z(r)-M(t)} < \sqrt{ 8
\alpha x } } \geq 1 - \frac{1}{8} . $$ %
Thus,
$$ \Pr \br{ \min_{0 \leq r \leq \k} \z(r) \geq M(t) - y } \geq
\Pr \br{ \min_{0 \leq r \leq x} \z(r) \geq M(t) - y } - \Pr \br{
\k > x } \geq \frac{1}{8} , $$ %
(where we have used the inequality $\Pr \br{ A \cap B } \geq \Pr
\br{A} - \Pr \br{ \ov{B}}$, valid for any events $A,B$).
\end{proof}

\begin{prop} \label{prop:  All particles stick without going below the layer M(t)-y} %
Let $G$ be a $d$-regular graph. Let $A_t$ be a \GCD cluster at
time $t$. Let $y \in \N$ and assume that
$$ \Pr \br{ \k(t+1) \leq y^2/4 } \geq \frac{1}{4} . $$
For $s > t$, define $H(s) = \z_s (\k(s))$; i.e. $H(s)$ is the
height of the layer at which the $s$-th particle sticks. Then, for
all $s>t$,
$$ \Pr \Br{ H(s) \geq M(t)- y }{ H(r) , \ t< r< s } \geq \frac{1}{8} , $$
for any values of $H(r)$, $t < r< s$.
\end{prop}

\begin{proof}
Let $s > t$.  Let $(g(\cdot),\z(\cdot)) =
(g_{s}(\cdot),\z_{s}(\cdot))$ be the walk the $s$-th particle
takes. That is, $g(0)$ is uniformly distributed in $G$, and $\z(0)
= M(s) \geq M(t)$. Set $k = \min \Set{ r \geq 0 }{ \z(r) = M(t)
}$, and set $k' = \min \Set{ r \geq k }{ (g(r),\z(r)) \in \partial
A_t }$.  $k$ is the first time the $s$-th particle is at the layer
$M(t)$ (this can be time $0$, e.g. if $M(s) = M(t)$).  $k'$ is the
first time after $k$ that the particle hits the outer boundary of
the cluster $A_t$.  Since $A_t \subseteq A_{s-1}$, we have that
$\k(s) \leq k'$. So,
\begin{eqnarray*}
    \Pr \Br{ H(s) \geq M(t)-y }{ A_{s-1} } & \geq & \Pr \Br{
    \min_{0 \leq r \leq \k(s)} \z(r) \geq M(t) -y }{ A_{s-1}
    } \\
    & \geq & \Pr \Br{ \min_{0 \leq r \leq k'} \z(r)
     \geq M(t) -y }{ A_{s-1} } \\
    & \geq & \Pr \Br{ \min_{0 \leq r \leq k'-k} \z(k+r) \geq M(t)-y }{ A_{s-1} } ,
\end{eqnarray*}
the last inequality following from the fact that for all $r< k$,
by definition, $\z(r) \geq M(t) \geq M(t) - y$.

Since the uniform distribution is the stationary distribution on
$G$, $g(k)$ is uniformly distributed in $G$. Thus, the walk
$(g(k+r),\z(k+r))$ has the same distribution as the walk
$(g_{t+1}(r),\z_{t+1}(r))$, and $k'-k$ has the same distribution
as $\k(t+1)$.  Using Proposition \ref{prop:  the particle sticks
without going below the layer M(t)-x} we now conclude
\begin{eqnarray*}
    \Pr \Br{ H(s) \geq M(t)-y }{ A_{s-1} } & \geq & \Pr \br{
    \min_{0 \leq r \leq \k(t+1)} \z_{t+1} (r) \geq M(t) -y }
    \geq \frac{1}{8} .
\end{eqnarray*}
Averaging over all $A \subset G \times \N$ such that $\Pr \Br{
A_{s-1} = A}{ H(r) , \ t<r<s } > 0$, we get that
$$ \Pr \Br { H(s) \geq M(t)- y }{ H(r) , \ t<r<s } \geq
\frac{1}{8} . $$
\end{proof}

\begin{prop} \label{prop:  if k(t+1) is small, then we build a wall higher than M-x} %
Let $G$ be a $d$-regular graph. Let $A_t$ be a \GCD cluster at
time $t$. Let $y \in \N$ and assume that
$$ \Pr \br{ \k(t+1) \leq y^2/4 } \geq \frac{1}{4} . $$ %

For $\ell \in \N$, define
$$ S_{\ell} = \min \Set{ s \geq t }{ L_s(\geq M(t)- y) = \ell +
L_t(\geq M(t)- y) } . $$ %
That is, $S_{\ell}$ is the first time that there are $\ell$ new
particles in the layers $\geq M(t)- y$ (so $S_0 = t$).

Then,
$$ \E \br{ S_{\ell} - t} \leq 8 \ell . $$
\end{prop}

\begin{proof}
The proof is similar to the proof of Proposition \ref{prop:
Expected time for ell particles above wall}.

By Proposition \ref{prop:  All particles stick without going below
the layer M(t)-y}, regardless of the previous particles, each
particle $s>t$ has probability at least $1/8$ to stick to a layer
$\geq M(t) - y$.  Thus, the expected time until there are $\ell$
particles above this layer is bounded by $8 \ell$.
\end{proof}

%
%

We now put everything together to prove Lemma \ref{lem:  if k(t+1)
is small then E[ T] leq n/loglog(n)}.  We show that if $\k(t+1)$
is small, then after a small amount of particles there is a high
layer with large load. Thus, after another small amount of
particles, the cluster reaches the new layer $M(t)$.

Set $M = M(t)$. Let
$$ T' = \min \Set{ s \geq t }{ \exists j \geq M-\mu \ : \ L_s(j)
\geq \frac{n}{\nu} } . $$ %
For $\ell \in \N$, define
$$ S_{\ell} = \min \Set{ s \geq t }{ L_s(\geq M-\mu) = \ell +
L_t(\geq M-\mu) } . $$ %
Consider the time $S_a$ for $a = \mu \lceil \frac{n}{\nu} \rceil$.
Consider the case where $T > S_a$. Then $M(S_a) = M$. At time
$S_a$, there are at least $a$ particles in the layers $\geq M-
\mu$, so
$$ a \leq \sum_{i=M-\mu}^{M(S_a)-1} L_{S_a}(i) = \sum_{i=M-\mu}^{M-1}
L_{S_a}(i) . $$ %
Thus, there exists $M-\mu \leq j \leq M-1$ such that $L_{S_a}(j)
\geq \frac{a}{\mu} \geq \frac{n}{\nu}$.  So $T' \leq S_a$.  We
conclude that if $T> S_a$ then $T' \leq S_a$.  In other words, we
have shown that $\min \set{T, T'} \leq S_a$. Hence, because it was
assumed that $\Pr \br{ \k(t+1) \leq \mu^2/4 } \geq \frac{1}{4}$,
by Proposition \ref{prop: if k(t+1) is small, then we build a wall
higher than M-x},
$$ \E \br{ \min \set{ T, T' } - t } \leq \E \br{ S_a -t } \leq 8a
\leq \frac{ 2n}{ \log \log (n) } + \frac{2 \log(n)}{\log \log (n)} . $$ %

Define the event
$$ B = \set{ \exists \ j \geq M(T') - \mu  \ : \ L_{T'}(j) \geq
\frac{n}{\nu} } . $$ %
By Lemma \ref{lem: Expected time to reach new layer, if there is a
wall}, we have that for large enough $n$,
$$ \E \br{ \1{B} \sr{ T - T' } } \leq
\frac{n}{4 \log \log n} .$$ %
We have that $T = \min \set{ T,T'} + \1{ T' < T } \sr{ T - T' }$.
Now, at time $T'$, we have a layer $j \geq M-\mu$ such that
$L_{T'}(j) \geq \frac{n}{\nu}$.  If $j < M(T') - \mu$ then $M(T')
> M$, and $T \leq T'$.  So the event $\set{ T' < T}$ implies the event $B$.
Thus, for large enough $n$,
$$ \E \br{ T -t} \leq \E \br{ \min \set{T,T'} - t} + \E \br{ \1{B} \sr{ T - T' } }
\leq \frac{5 n}{2 \log \log (n)} . $$
\end{proof}



\subsection{Particles take many steps}

In the previous section, we analyzed what happens when $\k(t+1)$
is ``small''.  This section is concerned with the case where
$\k(t+1)$ is ``large''.  The main goal of this section is proving
Lemma \ref{lem:  there are many steps, Pr of sticking to a new
layer} and Corollary \ref{cor:  when there are many steps, for RM
graphs}.  These are essential ingredients in the proof of Theorem
\ref{thm:  E [ T(M+1) - t ] leq n/loglog n}.

We begin with two technical lemmas:

\begin{lem} \label{lem:  G CDLA is distributed as G' CDLA}
Let $G$ be a graph.  Let $G'$ be the graph obtained from $G$ by
adding a self loop at each vertex.  That is,
$$ V(G') = V(G) \ , \ \quad \textrm{ and } \quad E(G') = E(G) \cup
\Set{ \set{v,v} }{ v \in G} . $$ %

Let $V = V(G \times \N) = V(G' \times \N)$.  Consider the \GCD and
$G'$-Cylinder-DLA processes.  Let $P_t(A,v) = \Pr \Br{ A_t = A
\cup \set{v} }{ A_{t-1} = A }$ where $\set{A_t}$ is a \GCD
process.  Let $P'_t(A,v) = \Pr \Br{ A_t = A \cup \set{v} }{
A_{t-1} = A }$ where $\set{A_t}$ is a $G'$-Cylinder-DLA process.

Then, for all $A \subseteq V$, $v \in V$, and all $t > 0$,
$$ P_t(A,v) = P'_t(A,v) . $$
\end{lem}

\begin{proof}
Assume that $A_t = A$.  We can couple the walk of the $(t+1)$-th
particle in both processes to hit the same vertex, as follows:

Denote by $L$ the set of self loops added to $G$ to form $G'$. Let
$\Set{ (g(r),\z(r)) }{ r \geq 0}$ be the walk of the $(t+1)$-th
particle, in the $G'$-Cylinder-DLA process. Define $\Gamma$ to be
the set of all $r > 0$ such that the step from $(g(r-1),\z(r-1))$
to $(g(r),\z(r))$ does not traverse one of the self loops in $L$.
For the \GCD process, let the $(t+1)$-th particle take the path
$\Set{ (g(r),\z(r)) }{ r \in \Gamma \cup \set{0} }$. This path has
the correct marginal distribution, as it is a simple random walk
on $G \times \N$. Note that both paths hit $\partial A_t$ at the
same vertex, since traversing a self loop does not move the
particle to a new vertex.
\end{proof}

\Rem If $G$ already has self loops, then by adding a self loop at
each vertex, we mean adding a new self loop, treated as different
from the original loop.  This only adds technical complications,
so we will not go into this issue.  The reader can treat all
graphs as not having self loops, though the results carry out to
graphs with self loops as well.

The important consequence of Lemma \ref{lem:  G CDLA is
distributed as G' CDLA} is that the Cylinder-DLA process does not
change if we let the particles perform a lazy random walk on $G
\times \N$.  This is needed to avoid technical complications that
arise from parity issues in bi-partite graphs.  The following
technical lemma is used to bypass this issue.

Recall our definition of the mixing time of a $d$-regular graph
$G$: Let $\set{g'_t}_{t \geq 0}$ be a lazy symmetric random walk
on $G$. The mixing time of $G$, is defined by
$$ \m = \m(G) \eqdef \min \Set{ t > 0 }{ \forall \ u,v \in G \ \forall \ s \geq t \ ,
\ \Pr \Br{ g'_s = u }{ g_0 = v } \geq \frac{1}{2 \abs{G}} } . $$ %

\begin{lem} \label{lem:  technical lemma for mixing time}
Let $G$ be a $d$-regular graph, and let $G'$ be the graph obtained
from $G$ by adding a self loop at each vertex, as in Lemma
\ref{lem: G CDLA is distributed as G' CDLA}. Let $\set{g_t}_{t
\geq 0}$ be a simple random walk on $G'$.  Then, for all $t \geq
\m(G)$, and all $u,v \in G$,
$$ \Pr \Br{ g_t = u }{ g_0 = v } \geq \frac{1}{2 \abs{G}} . $$
\end{lem}

\begin{proof}
This is immediate from the definition of $\m(G)$, and the fact
that $\set{g_t}_{t \geq 0}$ is distributed as a lazy symmetric
random walk on $G$.
\end{proof}

This completes the two technical lemmas we require.  Next we
introduce some notation.

Let $G$ be a graph.  Let $(g(0),\z(0)),(g(1),\z(1)),\ldots,$ be a
simple random walk on $G \times \N$.  For two times $r_1 < r_2$
denote
$$ r_1 \to r_2 \eqdef \set{ (g(r_1),\z(r_1)) , (g(r_1+1),\z(r_1+1)) ,
\ldots, (g(r_2),\z(r_2)) } . $$ %
$r_1 \to r_2$ is the path the walk takes between times $r_1$ and
$r_2$. Define
$$ L = \Set{ r > 0 }{ \z(r) = \z(0) } , $$
and assume that $L = \set{ \ell_1 < \ell_2 < \cdots }$. $L$ is the
set of times at which the walk visits the original layer. For $i
\geq 1$ define $\rho_i \eqdef \ell_{i-1} \to \ell_i$, where
$\ell_0 = 0$. We call $\rho_i$ an \emph{excursion}.   For $i \geq
1$ and $\alpha \in \R$, we say that $\rho_i = \ell_{i-1} \to
\ell_i$ is a \emph{positive $\alpha$-long excursion} if the
following conditions hold:
\begin{enumerate}
    \item $\z(\ell_{i-1}+1) = \z(0) + 1$; i.e. the excursion is on
    the positive side of the origin of the walk.

    \item The walk takes at least $\alpha$ steps in $G$ during the
    excursion; that is,
    $$ \sum_{r = \ell_{i-1}+1}^{\ell_i} \1{\z(r) = \z(r-1)} \geq
    \alpha . $$
\end{enumerate}
We stress that `$\alpha$-long' refers to the number of steps in
$G$, not the total length of the excursion.


%
%

\begin{lem} \label{lem:  there are many steps, Pr of sticking to a new
layer} %
Let $2 \leq d \in \N$.  There exist $c=c(d) > 0$ and $C=C(d) > 0$
such that for any $x \geq 1$ the following holds:

Let $G$ be a $d$-regular graph of size $\abs{G} = n$ and mixing
time $\m(G)$. Let $A_t$ be a \GCD cluster at time $t$. Recall that
$M(t)$ is the lowest empty layer at time $t$, and that $\k(t+1)$
is the number of steps the $(t+1)$-th particle takes before it
sticks.

Then,
$$ \Pr \Br{ M(t+1) > M(t) }{ A_t } > c \frac{x}{n \sqrt{\m(G)} }
\cdot \sr{ \frac{1}{2} - \Pr \br{ \k(t+1) \leq C x^2 } } .
$$
\end{lem}

\begin{proof} 

Let $G$ be a $d$-regular graph.  Let
$(g(0),\z(0)),(g(1),\z(1)),\ldots,$ be a simple random walk on $G
\times \N$. Let $\Set{ \rho_i = \ell_{i-1} \to \ell_i }{ i \geq
1}$ be the excursions of the walk.

First, we need to calculate the probability of a positive
$\alpha$-long excursion.

\begin{prop} \label{prop:  Pr of a positive long excursion}
For all $i \geq 1$ and any $2 \leq \alpha \in \R$, the probability
that $\rho_i$ is a positive $\alpha$-long excursion is greater
than $\frac{1}{12 (d+2) \sqrt{ \alpha } }$.
\end{prop}

\begin{proof}
Because of the Markov property, and the fact that $\z(\ell_i) =
\z(0)$ for all $i$, we get that $\Set{ \rho_i }{ i \geq 1}$ are
independent and identically distributed.  Thus, it suffices to
prove the proposition for $\rho = \rho_1$.

Fix $2 \leq \alpha \in \R$.  Set $m = \z(0)+1$.  So, the
probability that $\rho$ is a positive $\alpha$-long excursion is
equal to
\begin{eqnarray}
    & & \Pr \br{ \z(1) = m \ , \ \sum_{r=1}^{\ell_1} \1{ \z(r) =
    \z(r-1) } \geq \alpha } \nonumber \\
    & = & \Pr \br{ \z(1) = m } \Pr \Br{ \sum_{r=2}^{\ell_1} \1{
    \z(r) = \z(r-1) } \geq \alpha }{ \z(1) = m } ,
    \label{eqn:  Pr rho is a positive long excursion}
\end{eqnarray}
(we use the fact that $\z(1) = m \neq \z(0)$). Define
$$ \Gamma = \Set{ r > 1 }{ \z(r) = \z(r-1) } , \quad \textrm{ and } $$
$$ Z = \Set{ r > 1 }{ \z(r) \neq \z(r-1) } . $$
$\Gamma$ (respectively, $Z$) is the set of times at which the walk
moves in $G$ (respectively, $\N$). Let $g_1 = g(1)$ and let
$g_2,g_3,\ldots,$ be the walk $\Set{ g(r) }{ r \in \Gamma }$. So
$g_1,g_2\ldots,$ is distributed as a simple random walk on $G$,
starting at $g(1)$. Let $\z_1 = \z(1)$ and let $\z_2,\z_3,\ldots,$
be the walk $\Set{ \z(r) }{ r \in Z }$. So $\z_1,\z_2\ldots,$ is
distributed as a simple random walk on $\N$, starting at $\z(1)$.

Set $\alpha' = 8 \left\lceil \frac{\alpha}{2} \right\rceil$. For
$r > 1$, let $I_r = \1{\z(r) = \z(r-1)}$. Set
$$ \gamma = \sum_{r=2}^{\alpha'+1} I_r . $$
$\gamma$ is the sum of $\alpha'$ independent, identically
distributed Bernoulli random variables, with mean $\frac{d}{d+2}
\geq \frac{1}{2}$. Using the Chernoff bound (see e.g. Appendix A
in
\cite{ProbMethod}), %
\eqn{ \Pr \br{ \gamma < \frac{\alpha'}{4} } & \leq & 2 \exp \sr{ -
\frac{\alpha'}{8} } \leq \frac{2}{e} <
\frac{3}{4} . } %
$\gamma$ is independent of $\z(1)$, so
\begin{eqnarray} \label{eqn:  Pr [gamma geq alpha] , in excursion prop}
\Pr \Br{ \gamma \geq \alpha }{ \z(1) = m} & \geq & \Pr \Br{ \gamma
\geq \frac{\alpha'}{4} }{ \z(1) = m } > \frac{1}{4} .
\end{eqnarray}

Consider the walk $\z_1,\z_2, \ldots,\z_{\alpha'+1}$, conditioned
on the event $\z(1) = m$.  Define the event
$$ B = \set{ \z_2 \geq m \ , \ \z_3 \geq m \ , \ \ldots \ , \
\z_{\alpha'+1} \geq m } . $$ %
Conditioned on $\z(1) = m$, the walk
$\z_1,\z_2,\ldots,\z_{\alpha'+1}$ is a simple random walk on $\N$
starting at $\z_1 = m$. Using Corollary \ref{cor:  Pr of hitting
-1 after 2n steps},
$$ \Pr \Br{ B }{ \z(1) = m } \geq \Pr \Br{ \z_2 \geq m \ , \ \ldots \
, \ \z_{\alpha'+1} \geq m }{ \z(1) = m } = 2^{-\alpha'} { \alpha'
\choose \alpha' / 2} . $$ %
A careful application of Stirling's approximation gives
\begin{eqnarray} \label{eqn:  Pr [B | z(1) = m] ge , in excursion prop}
\Pr \Br{ B }{ \z(1) = m } & > & \frac{1}{3 \sqrt{ \alpha }} ,
\end{eqnarray}
for all $\alpha \geq 2$.

Since $\z(\ell_1) = \z(0) = m-1$, we have that, conditioned on
$\z(1) = m$, the event $B$ implies the event $\set{ \alpha'+1 \leq
\ell_1}$. Thus,
$$ \sum_{r=2}^{\ell_1} \1{ \z(r) = \z(r-1) } \geq \1{B}
\sum_{r=2}^{\alpha'+1} \1{ \z(r) = \z(r-1) } = \1{B} \gamma . $$
Now, $\gamma$ is independent of the event $B$, so, using
(\ref{eqn: Pr [gamma geq alpha] , in excursion prop})
and (\ref{eqn:  Pr [B | z(1) = m] ge , in excursion prop}), %
\eqn{ \Pr \Br{ \sum_{r=2}^{\ell_1} \1{
    \z(r) = \z(r-1) } \geq \alpha }{ \z(1) = m }
    & \geq & \Pr \Br{ B \ , \ \gamma \geq \alpha }{ \z(1) = m } \\
    & = & \Pr \Br{ B }{ \z(1) = m} \Pr \Br{ \gamma \geq \alpha }{
    \z(1) = m } \\
    & > & \frac{1}{12 \sqrt{ \alpha }} . } %
Plugging this into (\ref{eqn:  Pr rho is a positive long
excursion}), we have that the probability that $\rho$ is a
positive $\alpha$-long excursion is greater than $\frac{1}{12(d+2)
\sqrt{ \alpha } }$.
\end{proof}

The next proposition bounds from below the probability of sticking
to the layer $M(t)$ at each excursion.

\begin{prop} \label{prop:  good chance of sticking at each excursion} %
For all $i \geq 1$,
$$ \Pr \Br{ \k(t+1) = \ell_i }{ \k(t+1) > \ell_{i-1} } \geq \frac{c}{n
\sqrt{\m(G)} } , $$ %
where $c = c(d) > 0$ is a constant that depends only on $d$.
\end{prop}

\begin{proof}
It suffices to prove that for any $u,v \in G$
$$ \Pr \Br{ g(\ell_i) = u  }{ g(\ell_{i-1}) = v } \geq \frac{c}{n
\sqrt{ \m(G)} } , $$ %
for some constant $c=c(d)>0$, depending only on $d$.

Let $G'$ be the graph obtained from $G$ by adding a self loop at
each vertex. By Lemma \ref{lem:  G CDLA is distributed as G' CDLA}
we can assume that $(g(\cdot),\z(\cdot))$ is a walk on $G' \times
\N$.

Let $\Gamma = \Set{ \ell_{i-1} < r \leq \ell_i }{ \z(r) = \z(r-1)
}$, and let $\gamma = \abs{\Gamma}$.  Let $\Set{ h_r }{ 0 \leq r
\leq \gamma }$ be the walk $\Set{ g(\ell_{i-1}+r) }{ r \in \Gamma
\cup \set{0} }$.  $h_0,h_1,\ldots,h_\gamma$ is the walk measured
only when moving in the $G$-coordinate.  Note that conditioned on
$\Gamma$, the walk $h_0,h_1,\ldots,h_\gamma$ has the distribution
of a lazy random walk on $G$.   Thus by Lemma \ref{lem:  technical
lemma for mixing time}, we have that for any $u,v \in G$,
$$ \Pr \Br{ h_\gamma = u }{ \gamma \geq \m(G) \ , \
h_0 = v } \geq \frac{1}{2n} . $$ %
Note that if $\rho_i$ is a $\m(G)$-long excursion then $\gamma
\geq \m(G)$. Thus, for any $u,v \in G$, using Proposition
\ref{prop:  Pr of a positive long excursion},
\begin{eqnarray*}
    & & \Pr \Br{ g(\ell_i) = u }{ g(\ell_{i-1})  = v } \\
    & \geq & \Pr \Br{ h_\gamma = u }{ h_0 = v \ , \ \gamma \geq \m(G) }
    \Pr \Br{ \gamma \geq \m(G) }{ h_0 = v } \\
    & \geq & \frac{1}{2n} \cdot \frac{c}{ \sqrt{\m(G)} } .
\end{eqnarray*}
\end{proof}

Back to the proof of Lemma \ref{lem:  there are many steps, Pr of
sticking to a new layer}: Note that the events $\set{ \k(t+1) =
\ell_i }_{i=0}^{\infty}$ are pairwise disjoint, and that for every
$i \geq 0$, we have $\set{ \k(t+1) = \ell_i } \subset \set{ M(t+1)
> M(t) }$.  Thus, using Proposition \ref{prop:  good chance of
sticking at each excursion} we now have for any $x \geq 1$,
\begin{eqnarray} \label{eqn:  M(t+1) ge M(t) w.p. greater than
E[k]} %
    \Pr \br{ M(t+1) > M(t) } & \geq & \sum_{i=1}^{\infty} \Pr \br{
    \k(t+1) = \ell_i } \nonumber \\
    & = & \sum_{i=1}^{\infty}  \Pr \Br{ \k(t+1) =
    \ell_i }{ \k(t+1) > \ell_{i-1} } \Pr \br{ \k(t+1) > \ell_{i-1}
    } \nonumber \\
    & \geq & x \Pr \br{ \k(t+1) > \ell_x } \cdot \frac{c}{n
    \sqrt{\m(G)} } ,
\end{eqnarray}
for a constant $c=c(d)>0$ depending only on $d$.

Since for any $C > 0$,
$$ \Pr \br{ \k(t+1) > \ell_x } \geq \Pr \br{ \ell_x \leq C x^2 } -
\Pr \br{ \k(t+1) \leq C x^2 } , $$ %
we are left with proving that there exists $C>0$ such that $\Pr
\br{ \ell_x \leq C x^2 } \geq \frac{1}{2}$ for any $x \geq 1$.
Note that $\ell_x > C x^2$ implies that the number of times the
walk $\z(\cdot)$ visits the layer $M(t)$ up to time $\lceil C x^2
\rceil$ is less than $x$.  Thus by Lemma \ref{lem:  Lemma of zeros
of Lazy RW (new)}, there exists $C = C(d) > 0$ such that
$$ \Pr \br{ \ell_x > C x^2 } \leq \frac{1}{2}. $$
\end{proof}  

\begin{cor} \label{cor:  when there are many steps, for RM graphs} %
Let $2 \leq d \in \N$. There exist $n_0 = n_0(d)$ such that the
following holds for all $n > n_0$:

Let $G$ be a $d$-regular graph of size $n$, and mixing time
$$ \m(G) \leq \frac{\log^2(n)}{(\log \log (n))^5} . $$
Consider the \GCD process. Let $A_t$ be a \GCD cluster at time
$t$. Set
$$ \mu = \mu(n) = \left\lfloor \frac{\log(n)}{4 \log\log(n) } \right\rfloor .
$$ %
Assume that $\Pr \br{ \k(t+1) \leq \mu^2/4 } < \frac{1}{4}$. Then,
$$ \Pr \br{ M(t+1) > M(t) } > \frac{\log\log(n)}{n} . $$
\end{cor}

\begin{proof}
Let $C$ and $c$ be as in Lemma \ref{lem:  there are many steps, Pr
of sticking to a new layer}. We can choose $x \geq
\frac{\log(n)}{(\log\log(n))^{(3/2)} }$ such that $C x^2 \leq
\mu^2/4$ and $\frac{c x}{\sqrt{\m(G)}} \geq \log \log(n)$ for
large enough $n$.  Plugging this into Lemma \ref{lem:  there are
many steps, Pr of sticking to a new layer}, we get
$$ \Pr \br{ M(t+1) > M(t) } \geq \frac{\log \log(n)}{n} . $$
\end{proof}


\subsection{Proof of Theorem \ref{thm:  E [ T(M+1) - t ] leq  n/loglog
n}} \label{scn:  Proof of Theorem E [ T(M+1) - t ] leq n/loglog n} %

For convenience, we restate the Theorem:

\begin{thm*}[\ref{thm:  E [ T(M+1) - t ] leq  n/loglog n}]
Let $2 \leq d \in \N$. There exists $n_0 = n_0(d)$, such that the
following holds for all $n > n_0$:

Let $G$ be a $d$-regular graph of size $n$, and mixing time
$$ \m(G) \leq \frac{\log^2(n)}{(\log \log (n))^5} . $$
Let $A_t$ be a \GCD cluster at time $t$. Define
$$ T = \min \Set{ s > t }{ M(s) > M(t) } . $$
$T$ is the first time that a particle sticks to the empty layer,
$M(t)$. Then,
$$ \E \br{ T - t } \leq \frac{4 n}{\log \log n} . $$
\end{thm*}

\begin{proof}
Set $M = M(t)$ and
$$ \mu = \mu(n) = \left\lfloor \frac{ \log(n)}{4 \log \log (n)} \right\rfloor
, \quad \textrm{ and } \quad \nu = \nu(n) = \log(n) . $$

For $s \geq t$, define
$$ \alpha(s) = \Pr \Br{ \k(s+1) \leq \mu^2/4 }{ A_s } ,  $$
(which is random variable that is a function of $A_s$).  Define
$$ \tau = \min \Set{ s \geq t }{ \alpha(s) \geq \frac{1}{4} } . $$

Fix $s > t$, and $t+1 \leq r \leq s$. By Corollary \ref{cor:  when
there are many steps, for RM graphs}, there exists $n_0 = n_0(d)$
such that for all $n> n_0$,
\begin{eqnarray*}
    & & \Pr \Br{ M(r) = M(t) \ , \ \alpha(r) < \frac{1}{4} }{
    \forall \ t+1 \leq q \leq r-1 \ M(q) = M(t) \ , \ \alpha(q) <
    \frac{1}{4} } \\
    & \leq & \Pr \Br{ M(r) = M(t) }{
    \forall \ t+1 \leq q \leq r-1 \ M(q) = M(t) \ , \ \alpha(q) <
    \frac{1}{4} } \\
    & \leq & 1 - \frac{\log\log(n)}{n} .
\end{eqnarray*}
Thus, for all $s > t$,
\begin{eqnarray*}
    & & \Pr \br{ \min \set{ T,  \tau} > s } = \Pr \br{ T > s \
    , \ \tau > s } \\
    & \leq & \Pr \br{ \forall \ t+1 \leq r \leq s \ M(r) = M(t) \
    , \ \alpha(r) < \frac{1}{4} } \\
    & = & \prod_{r=t+1}^{s} \Pr \Br{ M(r) = M(t) \ , \ \alpha(r) < \frac{1}{4} }{
    \forall \ t+1 \leq q \leq r-1 \ M(q) = M(t) \ , \ \alpha(q) <
    \frac{1}{4} } \\
    & \leq & \sr{ 1 - \frac{\log\log(n)}{n} }^{s-t} .
\end{eqnarray*}
Since, $\Pr \br{ \min \set{ T, \tau} > t } \leq 1$, we get that
$$ \E \br{ \min \set{ T, \tau } - t } \leq \frac{n}{\log\log(n)}
. $$ %

Define
$$ T' = \left\{ \begin{array}{lr}
    \min \Set{ s > 0 }{ M(\tau+s) > M(\tau) } & \tau < \infty \\
    0 & \tau = \infty .
    \end{array} \right. $$
Then we have $T \leq \min \set{ T,\tau} + T'$. If $\tau = \infty$
then $\E \Br{ T' }{ \tau = \infty } = 0$. Assume that $\tau <
\infty$.  Then, at time $\tau$, we have that $\Pr \Br{ \k(\tau+1)
\leq \mu^2/4 }{ A_{\tau} } \geq \frac{1}{4}$. So, using Lemma
\ref{lem:  if k(t+1) is small then E[ T] leq n/loglog(n)},
$$ \E \Br{ T' }{ \tau < \infty } \leq \frac{5n}{2
\log\log(n)} , $$ %
and consequently,
$$ \E \br{ T' } < \frac{3n}{\log\log(n)} . $$
Thus, we conclude that
$$ \E \br{ T - t} \leq \E \br{ \min \set{ T,\tau} - t} + \E \br{
T' } < \frac{4n}{\log\log(n)} . $$
\end{proof}


%
%
%
%


\section{Density} \label{scn: Density}

\subsection{Definitions and Notation}

\begin{dfn}
Fix a graph $G$, and let $\set{A_t}$ be a \GCD process. Define the
\emph{cluster at infinity} by
$$ A_{\infty} = \bigcup_{t=0}^{\infty} A_t . $$

For $m \in \N$, define
$$ D(m) = \frac{1}{m n} \sum_{i=1}^m \abs{ A_{\infty} \cap G_i } . $$
$D(m)$ is the fractional amount of particles in the finite
cylinder $G \times \set{1,\ldots,m}$.

Define the \emph{density at infinity} by
\begin{equation} \label{eqn:  Density at infinity}
D = D_{\infty} = \lim_{m \to \infty} D(m) .
\end{equation}
\end{dfn}

Using standard arguments from ergodic theory it can be shown that
the limit in (\ref{eqn:  Density at infinity}) exists, and is
constant almost surely. Since $D(m)$ are bounded random variables,
we get by dominated convergence (see e.g. Chapter 9 in
\cite{Doob}):
$$ D = \E \br{ D} = \lim_{m \to \infty} \E \br{ D(m) } . $$


Recall the random times:
$$ T_m = \min \Set{ t \geq 0 }{ A_t \cap G_m \neq \emptyset } . $$
$T_m$ is the time the cluster first reaches the layer $m$.

\begin{thm} \label{thm:  E[D] leq lim (1/mn) E[T_m']}
Let $G$ be a $d$-regular graph of size $n$, and let $\set{A_t}$ be
a \GCD process. Let $D = D_{\infty}$ be the density at infinity,
and for all $m$ let
$$ T_m = \min \Set{ t \geq 0 }{ A_t \cap G_m \neq \emptyset } . $$

Then,
$$ D = \lim_{m \to \infty} \frac{1}{m n} \E \br{ T_{m} } . $$
\end{thm}


Theorem \ref{thm:  E[D] leq lim (1/mn) E[T_m']} relates the
density at infinity to the average growth rate. Theorem \ref{thm:
E[D] leq lim (1/mn) E[T_m']} is proved via the following
propositions.  The proof of the theorem is in Section \ref{scn:
Proof of theorem E[D] leq lim (1/mn) E[T_m']}.

\subsection{ }

The main objective of this section is Proposition \ref{prop:  Pr
[X(m)] leq p(m) 1/q }.  This proposition is the main observation
in proving Theorem \ref{thm:  E[D] leq lim (1/mn) E[T_m']}.

First we require some notation: For a \GCD process $\set{A_t}$,
recall $L_t(i) = \abs{ A_t \cap G_i }$, the load of the $i$-th
layer at time $t$. Define the load of the $i$-th layer at
infinity:
$$ L(i) = L_{\infty}(i) \eqdef \abs{ A_{\infty} \cap G_i } . $$

Define:
$$ L_t(\leq i) \eqdef \sum_{j=1}^i L_t(j) . $$ %
$$ L(\leq i) = L_{\infty}(\leq i) \eqdef \sum_{j=1}^i L(i) . $$ %
For $0 \leq t \leq \infty$, $L_t(\leq i)$ is the total load of all
layers below $i$, including $i$ but not including the $0$-layer.
(When indices become too small we write $L(t,\leq i)$ instead of
$L_t(\leq i)$.)

Also define $H(t) = \z_t(\k(t))$.  That is, $H(t)$ is the layer at
which the $t$-th particle sticks (the height of the $t$-th
particle).

The following proposition bounds the probability that a particle
sticks to a ``low'' layer.

\begin{prop} \label{prop:  prob of a particle to pass phi(m) is small} %
Fix $m < m' \in \N$. Let $G$ be a $d$-regular graph of size $n$,
with spectral gap $1-\lambda$ (i.e., $\lambda$ is the second
eigenvalue of the transition matrix of $G$). Consider the \GCD
process. Let $t
> T_{m'}$ and let $A_{t-1}$ be the \GCD cluster at time $t-1$.
Then,
$$ \Pr \Br{ H(t) \leq m }{ A_{t-1} }
< 3 \exp \sr{ - \frac{1-\lambda}{8n} (m'-m) } . $$
\end{prop}

\begin{proof}
Let $t > T_{m'}$. Let $M = M(t-1)$ and $\varphi = m'-m$.  Note
that
$$ M - 1 - m \geq M(T_{m'}) -1-m = m'-m = \varphi . $$
Let $(g(\cdot),\z(\cdot)) = (g_t(\cdot),\z_t(\cdot))$ be the walk
of the $t$-th particle. So $\z(0) = M$ and $g(0)$ is uniformly
distributed in $G$.  Let $k$ be the first step at which the walk
is at the layer $m$. That is, $k = \min \Set{ r > 0}{ \z(r) = m
}$. Let $\k = \k(t)$ be the step at which the particle sticks to
the cluster.

Note that the event $\set{ H(t) \leq m}$ implies the event $\set{
\k \geq k}$. Moreover, $\set{ \k \geq k }$ implies the event
$$ \set{ g(k-i) \not\in A_{t-1} \cap G_{\z(k-i)} \ , \ i = 1,2, \ldots, \varphi } . $$ %
Also, for all $1 \leq i \leq \varphi$ we have that $\abs{A_{t-1}
\cap G_{\z(k-i)} } \geq 1$ (because $\z(k-i) \leq m+i \leq M-1$).

Define
$$ \Gamma = \Set{ k-\varphi \leq r < k }{ \z(r) = \z(r-1) } , $$
and assume that
$$ \Gamma = \set{ r_1 < r_2 < \cdots < r_s }  $$
(note that $s =\abs{\Gamma}$ is a random variable). For $1 \leq i
\leq s$ let $g_i = g(r_i)$.  So $g_1,g_2,\ldots,g_s$ is
distributed as an $s$-step simple random walk on $G$, starting
from a uniformly chosen vertex.

For all $1 \leq i \leq s$ define $C_i = A_{t-1} \cap G_{\z(r_i)}$.
Thus, the event $\set{ H(t) \leq m }$ implies the event $\set{ g_i
\not\in C_i \ , \ i =1, 2, \ldots, s }$.

By Lemma \ref{lem:  RW stays in specific sets} we have that %
\eqn{ & & \Pr \Br{ g_i \not\in C_i \ , \ i =1, 2, \ldots, s } {
C_1,C_2, \ldots, C_s \ , \ s \geq \frac{1}{4} \varphi } \\
& \leq & \exp \sr{ - \frac{1-\lambda}{2n} \sum_{i=1}^s \abs{C_i} }
\leq \exp \sr{ - \frac{1-\lambda}{8n} \varphi } . } %
Hence, %
\eqn{ \Pr \Br{ H(t) \leq m }{ A_{t-1} } & \leq & \Pr \br{ g_i
\not\in C_i \ , \ i=1,2,\ldots, s } \\
& \leq & \Pr \br{ s < \frac{1}{4} \varphi } +
\exp \sr{ - \frac{1-\lambda}{8n} \varphi } . } %

Note that
$$ s = \sum_{i=1}^{\varphi} \1{ \z(k-i) = \z(k-i-1) } . $$
That is, $s$ is the sum of independent identically distributed
Bernoulli random variables, with mean $\frac{d}{d+2}\geq
\frac{1}{2}$.  Thus, using the Chernoff bound (see e.g. Appendix A
in \cite{ProbMethod}),
$$ \Pr \br{ s < \frac{1}{4} \varphi } < 2 \exp \sr{ -
\frac{\varphi}{8} } . $$ %
Thus, %
\eqn{ \Pr \Br{ H(t) \leq m }{ A_{t-1} } < 3 \exp \sr{ - \frac{1-\lambda}{8n} \varphi } . } %
\end{proof}

Consider the following event in the \GCD process:  Given a cluster
$A_t$, the next $\abs{G}$ particles appear in exactly the right
order so that they completely fill up the layer $M(t)$. (There is
always such an order; e.g. consider a spanning tree of $G$ rooted
at a vertex in $\partial A_t$.) Thus, an impassible ``wall'' is
created. Specifically, we are interested in the event that
$L_{t+n}(M(t)) = n$. The following proposition bounds from below
the probability of this event.

\begin{prop} \label{prop:  prob to get a wall immediately is large} %
Let $G$ be a $d$-regular graph of size $n$.  Consider the \GCD
process. Let $A_t$ be the \GCD cluster at time $t$. Then,
$$ \Pr \Br{ L_{t+n}(M(t)) = n }{ A_{t} } \geq
(d+2)^{-(n-1)} n^{-n} . $$
\end{prop}

\begin{proof}
Consider the following event $W$:  The $(t+1)$-th particle appears
at a vertex in $G_{M(t)}$ that is in $\partial A_t$.  Since there
is at least one such vertex, this happens with probability at
least $1/n$.  For $i=2,\ldots,n$, the $(t+i)$-th particle appears
at the layer $M(t)+1$, and moves to a vertex in $G_{M(t)}$ that is
in $\partial A_{t+i-1}$.  Since there is at least one such vertex,
the probability of this is at least $n^{-1} (d+2)^{-1}$, for each
$i=2,\ldots,n$.

Since the event $W$ implies that $L_{t+n}(M(t)) = n$, we have
$$ \Pr \Br{ L_{t+n}(M(t)) = n }{ A_{t} } \geq
(d+2)^{-(n-1)} n^{-n} . $$
\end{proof}

\begin{prop} \label{prop:  Pr [X(m)] leq p(m) 1/q }
Fix $m \in \N$. Let $\varphi = \varphi(m)$ be a positive integer,
and let $m' = m+\varphi$. Let $G$ be a $d$-regular graph of size
$n$, with spectral gap $1-\lambda$. Consider the \GCD process. Let
$X(m)$ be the event that there exists $t > T_{m'}$ such that $H(t)
\leq m$. That is, $X(m)$ is the event that a particle sticks to a
layer $\leq m$ after the cluster has reached the layer $m'$. Then,
$$ \Pr \br{ X(m) } \leq (d+2)^{n-1}(n+1) n^n \cdot
3 \exp \sr{ - \frac{1-\lambda}{8n} \varphi } . $$
\end{prop}

\begin{proof}
Fix $m \in \N$ and let $m' = m+\varphi(m)$.  Let $t > T_{m'}$. For
$i \in \N$ define the events
$$ W(t+i) = \set{ L_{t+i+n}(M(t+i)) = n } \ , \quad \textrm{ and } \quad
B(t+i) = \set{ H(t+i) \leq m } . $$ %
Set
$$ F(t+i) = \ov{B(t+i)} \cap \ov{W(t+i)} . $$

By Proposition \ref{prop:  prob of a particle to pass phi(m) is
small} we have that for all $i \geq n+1$,
$$ \Pr \Br{ B(t+i) }{ \forall \ 0 \leq j \leq i-(n+1) \ \ F(t+j) } \leq
3 \exp \sr{ - \frac{1-\lambda}{8n} \varphi(m) } . $$ %
By Proposition \ref{prop:  prob to get a wall immediately is
large} we have that for all $i \geq n+1$,
\begin{eqnarray*}
    & & \Pr \Br{ F(t+i) }{ \forall \ 0 \leq j \leq i-(n+1) \ \ F(t+j)
    } \\
    & \leq & \Pr \Br{ \ov{W(t+i)} }{ \forall \ 0 \leq j \leq i-(n+1) \ \ F(t+j)
    } \\
    & \leq & 1 - (d+2)^{-(n-1)} n^{-n} .
\end{eqnarray*}
Thus, for all $i \geq n+1$, %
\eqn{ & & \Pr \br{ B(t+i) \ , \ F(t+i-1) \ , \ \ldots \ , \ F(t) }
\\
& \leq & \Pr \Br{ B(t+i) }{ \forall \ 0 \leq j \leq i-(n+1) \ \
F(t+j) } \\
& & \times \prod_{\ell=1}^{\lfloor i/(n+1) \rfloor} \Pr \Br{
F(t+i- \ell(n+1)) }{ \forall \ 0 \leq j \leq i-(\ell+1)(n+1) \ \
F(t+j) } \\
& \leq & p(m) (1-q)^{\lfloor i/(n+1) \rfloor} , } %
where
$$ p(m) = 3 \exp \sr{ - \frac{1-\lambda}{8n} \varphi(m) } \quad
\textrm{ and } \quad q = (d+2)^{-(n-1)} n^{-n} . $$

Note that for all $t > T_{m'}$, the event $W(t)$ implies that
$\ov{B(t+i)}$ for all $i \geq 0$ (since the first $n$ particles
must stick to the layer $M(t) > m$, and after time $t+n$ no
particle can pass the layer $M(t) > m$). Thus, setting
$t=T_{m'}+1$, the event $X(m)$ implies that there exists $i \geq
0$ such that
$$ B(t+i) \cap \bigcap_{j=0}^{i-1} F(t+j) $$
occurs (i.e. take the first $i$ for which $B(t+i)$ occurs). So,
\begin{eqnarray*}
    \Pr \br{ X(m) } & \leq & \sum_{i=0}^{\infty}
    \Pr \br{ B(t+i) \ , \ F(t+i-1) \ , \ \ldots \ , \ F(t) } \\
    & \leq & \sum_{i=0}^{\infty} p(m) (1-q)^{\lfloor i/(n+1) \rfloor} \\
    & = & \sum_{\ell=0}^{\infty} (n+1) p(m) (1-q)^{\ell} = (n+1) p(m) \frac{1}{q} .
\end{eqnarray*}
\end{proof}

\subsection{Proof of Theorem \ref{thm:  E[D] leq lim (1/mn)
E[T_m']}} \label{scn:  Proof of theorem E[D] leq lim (1/mn)
E[T_m']} %

We restate the theorem:
\begin{thm*}[\ref{thm:  E[D] leq lim (1/mn) E[T_m']}]
Let $G$ be a $d$-regular graph of size $n$, and let $\set{A_t}$ be
a \GCD process. Let $D = D_{\infty}$ be the density at infinity,
and for all $m$ let
$$ T_m = \min \Set{ t \geq 0 }{ A_t \cap G_m \neq \emptyset } . $$

Then,
$$ D = \lim_{m \to \infty} \frac{1}{m n} \E \br{ T_{m} } . $$
\end{thm*}

\begin{proof}
Let $\varphi:\N \to \N$ be any function such that
$$ \lim_{m \to \infty} \varphi(m) = \infty \quad \textrm{ and }
\quad \lim_{m \to \infty} \frac{\varphi(m)}{m} = 0 . $$ %
For $m \in \N$ let $m' = m+\varphi(m)$.

Recall that $H(t) = \z_t(\k(t))$ is the height of the layer at
which the $t$-th particle sticks. For $m \in \N$ let $X(m)$ be the
event that there exists $t > T_{m'}$ such that $H(t) \leq m$.
Then, for all $\ell$, we have that
$$ \set{ L_{\infty}(\leq m) > \ell } \quad \textrm{ implies } \quad \set{
L(T_{m'},\leq m) > \ell } \cup X(m) . $$ %
This is because if $L(T_{m'},\leq m) \leq \ell$, then at least one
more particle is needed to stick at a layer $\leq m$ after time
$T_{m'}$, in order for $L_{\infty}(\leq m) > \ell$ to hold.

Thus, since $L_{\infty}(\leq m) \leq m n$, using Proposition
\ref{prop: Pr [X(m)] leq p(m) 1/q },
\begin{eqnarray*}
    \E \br{ L_{\infty}(\leq m) } & = & \sum_{\ell=0}^{\infty} \Pr \br{
    L_{\infty}(\leq m) > \ell } \\
    & = & \sum_{\ell=0}^{mn-1} \Pr \br{ L_{\infty}(\leq m) > \ell } \\
    & \leq & \sum_{\ell=0}^{\infty} \Pr \br{ L(T_{m'},\leq m) >
    \ell } + \sum_{\ell=0}^{mn-1} \Pr \br{ X(m) } \\
    & = & \E \br{ L(T_{m'},\leq m) } + m n \cdot p(m) \frac{n+1}{q} ,
\end{eqnarray*}
for
$$ p(m) = 3 \exp \sr{ - \frac{1-\lambda}{8n} \varphi(m) } \quad
\textrm{ and } \quad q = (d+2)^{-(n-1)} n^{-n} . $$

Note that $L_t(\leq m) \leq t$ for all $t$, so
$$ \E \br{ L(T_{m'}, \leq m) } \leq \E \br{ T_{m'} } . $$
Also,
$$ \lim_{m \to \infty} p(m) \frac{n+1}{q} = 0 . $$
So,
\begin{eqnarray*}
    \E \br{ D } & = & \lim_{m \to \infty} \E \br{ D(m) } = \lim_{m
    \to \infty} \frac{1}{m n} \E \br{ L_{\infty}(\leq m) } \\
    & \leq & \lim_{m \to \infty} \frac{1}{m n} \E \br{ T_{m'} } +
    \lim_{m \to \infty} p(m) \frac{n+1}{q} \\
    & = & \lim_{m \to \infty} \frac{1}{m n} \E \br{ T_{m'} } .
\end{eqnarray*}
Since for all $k' > k$, $\E \br{ T_{k'} - T_k } \leq n(k'-k)$, we
have that
$$ \E \br{ T_{m'} } = \E \br{ T_m } + \E \br{ T_{m'} - T_m } \leq
\E \br{T_m} + \varphi(m) n . $$ %
Thus,
\begin{equation} \label{eqn:  upper bound E[D]}
D = \E \br{ D } \leq \lim_{m \to \infty} \frac{1}{mn} \E \br{ T_m
} + \lim_{m \to \infty} \frac{\varphi(m)}{m} = \lim_{m \to \infty}
\frac{1}{mn} \E \br{ T_m } .
\end{equation}

Note that for all $m$,
$$ T_m = L(T_m, \leq m) \leq L_{\infty}(\leq m) , $$
so $\E \br{ T_m } \leq \E \br{ L_{\infty}(\leq m)}$.  Thus,
\begin{equation} \label{eqn:  lower bound E[D]}
\lim_{m \to \infty} \frac{1}{m n} \E \br{ T_m} \leq \lim_{m \to
\infty} \frac{1}{m n} \E \br{ L_{\infty}(\leq m)} = \E \br{D} = D
.
\end{equation}
(\ref{eqn:  upper bound E[D]}) and (\ref{eqn:  lower bound E[D]})
together give equality.
\end{proof}

\subsection{Density of Cylinder-DLA with transitive base}

In this section we assume that $G$ is vertex transitive; i.e. for
any $u,v \in G$ there exists an automorphism (of graphs)
$\varphi_{uv}:G \to G$ such that $\varphi(u) = v$.

\begin{thm} \label{thm:  for any v.t. graph D leq 2/3}
Let $G$ be a vertex transitive graph.  Let $\set{A_t}$ be the \GCD
process. Let $D = D_{\infty}$ be the density at infinity. Then,
$$ D \leq \frac{2}{3} . $$
\end{thm}

The key to proving Theorem \ref{thm:  for any v.t. graph D leq
2/3} is Lemma \ref{lem:  Pr [ X(t-1) in S(t) ] geq 2d+1/d+1}
below.  The proof of the Theorem follows the proof of the Lemma.

\begin{lem} \label{lem:  Pr [ X(t-1) in S(t) ] geq 2d+1/d+1}
Let $G$ be a vertex transitive graph.  Let $A_{t-1}$ be a \GCD
cluster at time $t-1$.  Then,
$$ \Pr \Br{ M(t) > M(t-1) }{ A_{t-1} } \geq \frac{2d+2}{(d+2)n} . $$
\end{lem}

\begin{proof}
Recall $(g_t(\cdot),\z_t(\cdot))$ is the walk of the $t$-th
particle, so $g_t(0)$ is uniformly distributed in $G$, and
$\z_t(0) = M(t-1)$.

%

Define $\X(t-1)$ to be the newest particle in the top layer of the
cluster $A_{t-1}$.  That is, if $A_{t-1} \cap G_{(M(t-1)-1)} =
\set{v_1,\ldots, v_\ell}$, then $\X(t-1)$ is the vertex $v_i$ that
is the last vertex to join the cluster.

Note that $(\X(t-1),M(t-1)) \in \partial A_{t-1}$.

Because the graph $G$ is vertex transitive, we get that $\X(t-1)$
is uniformly distributed in $G$.  Moreover, $\X(t-1)$ depends only
on the clusters $A_{t-1}, \ldots, A_0$, and is independent of the
walk $(g_t(\cdot),\z_t(\cdot))$.

Define $S(t)$ to be the set of vertices in $G$ that the walk
$(g_t(\cdot),\z_t(\cdot))$ visits before leaving the layer
$M(t-1)$.  That is:
$$ \tau = \min \Set{ r > 0 }{ \z_t(r) \neq M(t-1) } ,$$
$\tau$ is the first step the $t$-th particle is not in the layer
$M(t-1)$.
$$ S(t) = \Set{ v \in G }{ \exists \ 0 \leq r \leq \tau-1 \ : \
g_t(r) = v } . $$

\Claim For all $t>1$,
$$ \Pr \br{ \X(t-1) \in S(t) } = \frac{1}{n} \E \br{ \abs{S(t)} } . $$

\begin{proof}
For any $u \in G$,
$$ \Pr \Br{ \X(t-1) = u }{ S(t) } = \frac{1}{n} . $$
Consequently,
\begin{eqnarray*}
    \Pr \br{ \X(t-1) \in S(t) } & = & \sum_S \Pr \Br{ \X(t-1) \in S }{
    S(t) = S } \Pr \br{ S(t) = S } \\
    & = & \sum_S \sum_{v \in S} \Pr \Br{ \X(t-1) = v }{ S(t) = S }
    \Pr \br{ S(t) = S } \\
    & = & \frac{1}{n} \sum_S \abs{S} \Pr \br{S(t) = S } = \frac{1}{n} \E \br{ \abs{ S(t) } } .
\end{eqnarray*}
\end{proof}

Recall that $d$ is the degree of $G$.

\Claim For all $t > 0$,
$$ \E \br{ \abs{S(t)} } \geq \frac{2d+2}{d+2} . $$

\begin{proof}
Let $R(k)$ denote the range of a $k$-step random walk on $G$.
Then, for $s \geq 2$,
\begin{eqnarray*}
    \Pr \br{ \abs{ S(t)} = s } & = & \sum_{k=s-1}^{\infty} \Pr
    \br{ R(k) = s } \sr{ \frac{d}{d+2} }^k \frac{2}{d+2} .
\end{eqnarray*}
For $s = 1$,
$$ \Pr \br{ \abs{ S(t) } =1 } = \frac{2}{d+2} . $$
Thus,
\begin{eqnarray} \label{eqn:  E [ |S(t)| ] geq sum E [R(k)]}
    \E \br{ \abs{S(t)} } & = & \frac{2}{d+2} + \sum_{s \geq 2} \sum_{k \geq s-1}
    \frac{2}{d+2} \sr{ \frac{d}{d+2} }^k s \Pr \br{ R(k) = s } \nonumber \\
    & = & \frac{2}{d+2} + \frac{2}{d+2} \cdot \sum_{k \geq 1}
    \sum_{s=2}^{k+1} s \Pr \br{ R(k) = s } \sr{ \frac{d}{d+2} }^k \nonumber \\
    & = & \frac{2}{d+2} \cdot \sum_{k \geq 0} \E
    \br{ R(k) } \sr{ \frac{d}{d+2} }^k
\end{eqnarray}
Substitute in (\ref{eqn:  E [ |S(t)| ] geq sum E [R(k)]}) the
naive bound $R(k) \geq 2$ for all $k \geq 1$:
\begin{eqnarray*}
    \E \br{ \abs{S(t)} } & \geq & \frac{2}{d+2} \cdot \br{ 1 + 2
    \sum_{k=1}^{\infty} \sr{ \frac{d}{d+2} }^k } \\
    & = & \frac{2}{d+2} \cdot \br{ 1 + 2 \frac{d}{d+2} \frac{d+2}{2} } =
    \frac{2d+2}{d+2} .
\end{eqnarray*}
\end{proof}

Thus, using the claim, we have that for all $t>1$,
$$ \Pr \br{ \X(t-1) \in S(t) } \geq \frac{2d+2}{(d+2)n} . $$
The lemma now follows from the fact that the event $\set{ \X(t-1)
\in S(t) }$ implies $\set{ S(t) \cap \partial A_{t-1} \neq
\emptyset}$.  So $\set{ \X(t-1) \in S(t) }$ implies that the
$t$-th particle sticks to the layer $M(t-1)$.
\end{proof}

\begin{proof}[Proof of Theorem \ref{thm:  for any v.t. graph D leq
2/3}] %
For $m \in \N$ recall
$$ T_m = \min \Set{ t \geq 0 }{ A_t \cap G_m \neq \emptyset } . $$
By Lemma \ref{lem:  Pr [ X(t-1) in S(t) ] geq 2d+1/d+1}, for all
$m$,
$$ \E \br{ T_m } \leq m \frac{(d+2)n}{2d+2} . $$
Thus,
$$ \frac{1}{mn} \E \br{ T_m } \leq \frac{d+2}{2d+2} . $$
Plugging this into Theorem \ref{thm:  E[D] leq lim (1/mn)
E[T_m']}, we have
$$ D \leq \frac{d+2}{2d+2} \leq \frac{2}{3} . $$
\end{proof}

\subsection{Density of Cylinder-DLA with quickly mixing base}

In this section we combine two main results:  For a family of
graphs $\set{G_n}$ with small mixing time, we show that since the
\GCD process grows arms, the densities at infinity tend to $0$ as
$n$ tends to infinity.  Formally:

\begin{thm} \label{thm: D(n) leq 3/log log n for RM family of
graphs} %
Let $2 \leq d \in \N$.  Let $\set{G_n}$ be a family of $d$-regular
graphs such that $\lim_{n \to \infty} \abs{G_n} = \infty$, and for
all $n$,
$$ \m(G_n) \leq \frac{\log^2 \abs{G_n} }{(\log \log \abs{G_n})^5} . $$
For all $n$ let $D(n)$ be the density at infinity of the
$G_n$-Cylinder-DLA process. Then,
$$ \lim_{n \to \infty} D(n) = 0 . $$
\end{thm}

\begin{proof}
There exists $n_0 = n_0 (d)$ such that the following holds for all
$n>n_0$:

Set $G = G_n$ and consider $\set{A_t}$, a \GCD process. By Theorem
\ref{thm:  E [ T_m ] leq o(nm)}, for all $m$,
$$ \E \br{ T_m } < m n \frac{4}{\log \log n} . $$
Thus, using Theorem \ref{thm:  E[D] leq lim (1/mn) E[T_m']},
$$ D(n) \leq \frac{4}{\log \log n} , $$
for all $n > n_0$. Thus,
$$ \lim_{n \to \infty} D(n) = 0 . $$
\end{proof}

\section{Lower bound on the growth rate} \label{scn:  lower bound section}

In this section we prove a lower bound on the expected growth rate
of the \GCD cluster, provided that the spectral gap is at least
$\abs{G}^{-2/3}$. This regime of the spectral gap includes graphs
with small mixing time as in Theorem \ref{thm:  E [ T_m ] leq
o(nm)}, and many more natural families of graphs such as discrete
cubes and tori of dimension at least $3$.

\begin{thm} \label{thm:  lower bound on speed}
Let $2 \leq d \in \N$. There exists $n_0 = n_0(d)$, such that the
following holds for all $n > n_0$:

Let $G$ be a $d$-regular graph such that
$$ \abs{G} = n \quad \textrm{ and } \quad 1-\lambda \geq n^{-2/3} , $$
where $1-\lambda$ is the spectral gap of $G$. Consider
$\set{A_t}$, a \GCD process. For $m \in \N$, define
$$ T_m = \min \Set{ t \geq 0 }{ A_t \cap G_m \neq \emptyset } . $$
$T_m$ is the time the cluster first reaches the layer $m$.

Then, for all $m$,
$$ \E \br{ T_m} > C m n^{1/20} , $$
where $C$ is some constant that depends only on $d$.
\end{thm}

\begin{proof}
Fix $t>0$, and let $A_{t-1}$ be the \GCD cluster at time $t-1$.

\Claim There exists a constant $C = C(d)$ (that depends on $d$)
such that for all $t > 0$,
$$ \Pr \Br{ M(t) > M(t-1) }{ A_{t-1} } <
\frac{C \abs{\partial A_{t-1} \cap G_{M(t-1)} } }{n^{1/10}} . $$ %

\begin{proof}
Let $(g(\cdot),\z(\cdot)) = (g_t(\cdot),\z_t(\cdot))$ be the walk
of the $t$-th particle. Set
$$ L = \Set{ r > 0 }{ \z(r) = \z(0) } = \set{ \ell_1 < \ell_2 <
\cdots } , $$ %
and let $\rho_i = \ell_{i-1} \to \ell_i$ be the excursions of the
walk.  For $2 \leq \alpha \in \R$, let $p(\alpha)$ be the
probability that an excursion is a negative $\alpha$-long
excursion; that is $p(\alpha)$ is the probability that
$$ \z(\ell_{i-1}+1) = \z(0) - 1 \quad \textrm{ and } \quad \sum_{r = \ell_{i-1}+1}^{\ell_i}
\1{\z(r) = \z(r-1)} \geq \alpha . $$ %
(This is independent of $i$.) By symmetry and Proposition
\ref{prop: Pr of a positive long excursion}, we have that
$p(\alpha) > (1/c(d)) \alpha^{-1/2}$, where $c(d) = 12(d+2)$.

Fix $2 \leq \alpha \in \R$, and set $p = p(\alpha)$.

For an integer $k \in \N$, let $N(k)$ denote the number of
negative $\alpha$-long excursions out of the first $k$ excursions.
So $N(k) = \sum_{i=1}^k I_i(\alpha)$, where $I_i(\alpha)$ is the
indicator of the event that $\rho_i$ is a negative $\alpha$-long
excursion. Since $\set{I_i(\alpha)}$ are indpendent, we have by
Chebychev's inequality that
$$ \Pr \br{ N(k) \leq \frac{1}{2} p k } \leq \frac{4}{pk} < 4 c(d) \frac{\sqrt{\alpha}}{k} . $$
Let $Z = Z(k)$ be the number of times up to $\ell_k$ the walk
moves in $G$ while on the negative side of $\z(0)$; i.e.,
$$ Z(k) = \sum_{r=1}^{\ell_k} \1{ \z(r) = \z(r-1) } \1{ \z(r) < \z(0) } . $$
We have that $Z \geq \alpha \cdot N(k)$ (since each negative
$\alpha$-long excursion contributes at least $\alpha$ to the sum).
Thus,
$$ \Pr \br{ Z \leq \frac{ \sqrt{\alpha} k}{2 c(d)} } \leq
\Pr \br{ Z \leq \frac{\alpha}{2} p k } < 4 c(d) \frac{\sqrt{\alpha}}{k} . $$


Set $A = \partial A_{t-1} \cap G_{M(t-1)}$. Set $B = \partial
A_{t-1} \setminus A$. That is, $B$ is the set $\partial A_{t-1}$
with the highest layer removed.

For all $r \geq 0$ let $C_r = B \cap G_{\z(r)}$. Note that if
$\z(r) < \z(0)$ then $\abs{C_r} \geq 1$ (because any layer below
$M(t-1)$ contains at least one particle). Define a simple random
walk on $G$ by $h_0 = g(0)$ and
$$ \set{ h_1, h_2, \ldots, h_s } = \Set{ g(r) }{ \z(r) = \z(r-1)
\ , \ 1 \leq r \leq \ell_k } . $$ %

Let $F = F(\ell_k)$ be the event that the particle does not hit
the set $B$ up to time $\ell_k$.  That is,
$$ F = \set{ \forall 0 \leq r \leq \ell_k \ \ (g(r),\z(r)) \not\in B
} . $$ %
Conditioned on a specific path $\z(0), \z(1), \ldots, \z(\ell_k)$,
and on $A_{t-1}$, we have that $h_0,h_1,\ldots,h_s$ is distributed
as a simple random walk on $G$. Using Lemma \ref{lem:  RW stays in
specific sets}, that
\begin{eqnarray*}
    & & \Pr \Br{ F }{ \z(0), \ldots, \z(\ell_k) \ , \ A_{t-1} } \\
    & \leq & \Pr \Br{ \forall \ 0 \leq r \leq s \ h_r \not\in C_r
    }{ \z(0), \ldots, \z(\ell_k) } \\
    & \leq & \exp \sr{-  \frac{1-\lambda}{2 n} \sum_{r=1}^s C_r }
    \leq \exp \sr{ - \frac{1-\lambda}{2 n} Z } .
\end{eqnarray*}
Thus, averaging over all possible paths $\z(0), \z(1), \ldots,
\z(\ell_k)$, we have that
\begin{eqnarray*}
\Pr \Br{ F }{ A_{t-1} } & < & \Pr \br{ Z \leq \frac{\alpha}{2} p k
} + \exp \sr{ - \frac{1-\lambda}{4 n} \alpha p k  } \\
& < & 4 c(d) \frac{\sqrt{\alpha}}{k} + \exp \sr{ -
\frac{1-\lambda}{4 c(d) n} \sqrt{\alpha} k } .
\end{eqnarray*}
Note that the event $\set{ \k(t) > \ell_k }$ implies the event
$F$, so we have that
$$ \Pr \Br{ \k(t) > \ell_k }{ A_{t-1} } < 4 c(d) \frac{\sqrt{\alpha}}{k} + \exp \sr{ -
\frac{1-\lambda}{4 c(d) n} \sqrt{\alpha} k } . $$ %

On the other hand, consider the times $\ell_0, \ell_1, \ldots,
\ell_k$. Since $\partial A_{t-1} \cap G_{\z(0)} = \partial A_{t-1}
\cap G_{M(t-1)} = A$, we have by a union bound,
\begin{eqnarray*}
\Pr \Br{ \exists \ x \in A \ , \ \exists \ 0 \leq i \leq k \ : \
(g(\ell_i),\z(\ell_i)) = x }{ A_{t-1} } & \leq & \frac{\abs{A}
(k+1)}{n} .
\end{eqnarray*}

Now, the event $\set{M(t) > M(t-1)}$ implies that there exists $i
\geq 0$ such that the particle does not stick to $\partial
A_{t-1}$ before time $\ell_i$, and $(g(\ell_i),\z(\ell_i)) = x$
for some $x \in A$. Thus, we have for all $2 \leq \alpha \in \R$
and all $k \in \N$,
\begin{eqnarray*}
& & \Pr \Br{ M(t) > M(t-1) }{ A_{t-1} } \\
& \leq & \Pr \Br{ \k(t)
> \ell_k }{ A_{t-1} } + \Pr \Br{ \exists \ x \in A \ , \ \exists \ 0 \leq i \leq k  \ : \
(g(\ell_i),\z(\ell_i)) = x }{ A_{t-1} } \\
& < & 4 c(d) \frac{\sqrt{\alpha}}{k} + \exp \sr{ -
\frac{1-\lambda}{4 c(d) n} \sqrt{\alpha} k } + \frac{\abs{A}
(k+1)}{n} .
\end{eqnarray*}


Set $\eps = 1/10$, $k = n^{1-\eps}$, $\alpha = n^{2-4\eps}$. Then,
if $1-\lambda \geq \frac{1}{n^{2/3}}$, we have that for large
enough $n$ (depending on $d$),
$$ \Pr \Br{ M(t) > M(t-1) }{ A_{t-1} } < \frac{C \abs{A} }{n^{\eps}} , $$ %
for some constant $C = C(d)$.
\end{proof}

Back to the proof of Theorem \ref{thm:  lower bound on speed}: Fix
$m > 0$, and consider the time $T_m$.  Note that for all $1 \leq j
\leq n$,
$$ \abs{ \partial A_{T_m+j-1} \cap G_{M(T_m+j-1)} } \leq j , $$
(because at most $j$ particles could have stuck to the layer
$M(T_m+j-1)-1$ by time $T_m+j-1$).  Thus, for all $1 \leq j \leq
n$ we have that for $C$ and $\eps$ as above
$$ \Pr \Br{ M(T_m+j) = M(T_m+j-1) }{ A_{T_m+j-1} } > 1 - \frac{C
j}{n^\eps} . $$ %
This implies that for $\lambda < \frac{n^{\eps}}{C}$, %
\eqn{ \Pr \br{ T_{m+1} - T_m > \lambda } & > &
\prod_{j=1}^{\lambda} \sr{ 1 - \frac{C j}{n^{\eps}} } \geq \sr{ 1
- \frac{C \lambda}{n^{\eps}} }^\lambda , } %
and so, there exists a constant $C'$ (depending on $C$) such that
for $\lambda = \lceil n^{\eps/2} \rceil$, %
\eqn{ \E \br{ T_{m+1} - T_m } & > & \lambda \Pr \br{ T_{m+1} - T_m
> \lambda } > C' n^{\eps/2} . } %
Hence, we get that for all $m \geq 2$,
$$ \E \br{ T_m } > \frac{C'}{2} m n^{\eps/2} . $$
\end{proof}

For completeness, we state the immediate Corollary of Theorems
\ref{thm:  lower bound on speed} and \ref{thm:  E[D] leq lim
(1/mn) E[T_m']}.

\begin{cor}
Let $2 \leq d \in \N$. There exists $n_0 = n_0(d)$, such that the
following holds for all $n > n_0$:

Let $G$ be a $d$-regular graph such that
$$ \abs{G} = n \quad \textrm{ and } \quad 1-\lambda \geq n^{-2/3} , $$
where $1-\lambda$ is the spectral gap of $G$. Consider
$\set{A_t}$, a \GCD process. Let $D_{\infty}$ be the density at
infinity.  Then, for some constant $C$ that depends only on $d$,
$$ D_{\infty} \geq \frac{C}{n^{19/20} } . $$
\end{cor}


\section{Further research directions}

The results and methods in this paper raise a few natural
questions:

\begin{enumerate}

\item Let $G$ be a $d$-regular graph.  Let $H$ be obtained from
$G$ by only adding edges to $G$, so that $V(H) = V(G)$ and $H$ is
$(d+1)$-regular.  Is there monotonicity in the expected speed of
the cluster on the Cylinder-DLA processes with base $G$ and with
base $H$.  That is, let $T^G_m$, respectively $T^H_m$, be the
first time the cluster reaches the layer $m$ in the \GCD,
respectively $H$-Cylinder-DLA, process.  Is it true that $\E \br{
T^G_m } \geq \E \br{ T^H_m }$ for all $m$?

\item Consider a \GCD process, started with $A_0 = \set{x_0}$ for
a specific vertex $x_0 \in G$.  Let $\tau$ be the mixing time of a
simple random walk on $G$ (i.e. the time it takes for a simple
random walk to come close in total-variation distance to the
stationary distribution). For $m > 0$, let $x_m \in G$ be the
vertex in $G$ that is the first vertex in the layer $m$ that a
particle sticks to. In our notation above $x_m = v$ such that
$A_{T_m} \cap G_m = \set{(v,m)}$.  How long does it take for the
distribution of $x_m$ to be close to the uniform distribution?
Does there exist a constant $c$ such that $x_{c \tau}$ is close to
being uniformly distributed on $G$?

\item Directed \GCD : Consider a model of \GCD where particles
cannot move to layers above, only to layers below or in their
current layer.  Is the density of directed \GCD always greater
than undirected?  Are there graphs $G$ for which these quantities
are of the same order?  Are there graphs for which the ratio
between the density of undirected \GCD and directed \GCD goes to
$0$ as the size of $G$ goes to infinity?

The model of directed \GCD can also be generalized to a model
where particles move up with probability $\frac{\alpha}{d+2}$ and
down with probability $\frac{2-\alpha}{d+2}$ (and to a neighbor in
the current layer with probability $\frac{1}{d+2}$), for some
$\alpha < 1$.  Thus, there is a drift down.  The same questions
can be asked of this model.

We remark that some of our results still hold in directed \GCD.
Mainly, Lemma \ref{lem:  if k(t+1) is small then E[ T] leq
n/loglog(n)} (that states that if the particle takes a small
amount of steps to stick, then the expected time to reach the new
layer is small,) still holds with the assumption that $\Pr \br{
\k(t+1) \leq \mu } \geq \frac{1}{4}$.

\item The \GCD process, is of course not a stationary process
(since $A_t \subset A_{t+1}$ for all $t$).  But, each time a
``wall'' is built (i.e. $L_{t+n}(M(t)) = n$, see Proposition
\ref{prop:  prob to get a wall immediately is large}), we start
the cluster again, independently of the cluster below the wall. If
we identify clusters that are the same above walls, we get a
stationary Markov chain on clusters.  Our analysis throughout this
paper in some sense evades this stationary distribution.  It would
be interesting if some properties of the cluster generated under
the stationary distribution could be worked out.  Perhaps,
calculating properties of the ``typical cluster'' could help
improve the results of this paper (e.g., reduce the spectral gap
required to grow arms).

\item As stated in the introduction, DLA on a cylinder suggests
studying the problem of ``clogging''.  That is, run a \GCD process
for some graph $G$.  Let $T$ be the (random) time at which the
cluster clogs the cylinder.  That is, $T$ is the first time at
which there exists a layer such that no particle can pass this
layer; i.e.,
$$ T = \min \Set{ t > 0 }{ \exists \ m > 1 \ : \ \Pr \br{ H(t)
\leq m } = 0 } . $$ Provide bounds on $\E \br{ T }$.  How is $T$
distributed?

\end{enumerate}


\appendix

\section{Random Walks on $\Z$}

We collect some facts about a simple random walk on $\Z$, $S(n)$,
starting at $S(0)= 0$.

The following is Theorem 9.1 of \cite{Revesz}:
\begin{lem}
Let
$$ \rho(1) = \min \Set{ i \geq 0 }{ S(i) = 1 } . $$
Then, for all $n$,
$$ \Pr \br{ \rho(1) > 2n } = 2^{-2n} { 2n \choose n } . $$
\end{lem}

\begin{cor} \label{cor:  Pr of hitting -1 after 2n steps}
For all $n$,
$$ \Pr \br{ \forall \ 1 \leq i \leq 2n \ , \ S(i) \geq 0 } = 2^{-2n}
{2n \choose n} . $$
\end{cor}

\begin{proof}
Let
$$ \tau = \min \Set{ i \geq 0 }{ S(i) = -1 } . $$
By symmetry, $\tau$ has the same distribution as $\rho(1)$ above.
Thus, for all $n$,
$$ \Pr \br{ \forall \ 1 \leq i \leq 2n \ , \ S(i) \geq 0 } =
   \Pr \br{ \tau > 2n } = \Pr \br{ \rho(1) > 2n } = 2^{-2n} {2n
   \choose n} . $$
\end{proof}

The following is Theorem 9.3 of \cite{Revesz}:
\begin{lem} \label{lem:  Pr[ L(2n) < m] }
Let $L(n)$ be the number of times the walk has visited $0$, i.e.
$$ L(n) = \abs{ \Set{ 1 \leq i \leq n }{ S(i) = 0 } } . $$
Then for $m \leq n$,
\begin{equation} \label{eqn:  Pr L(2n) < m }
\Pr \br{ L(2n) < m } = 2^{-2n} \sum_{j=0}^{m-1} 2^j { 2n-j \choose
n } .
\end{equation}
\end{lem}

\begin{cor} \label{cor: Stirling on L(n)}
For $L(n)$ as above, and $m \leq n/2$,
$$ \Pr \br{ L(n) < m } < \frac{m}{ \sqrt{n-2m} } . $$
\end{cor}

\begin{proof}
This is a careful application of Stirling's approximation to
(\ref{eqn:  Pr L(2n) < m }).
\end{proof}

\begin{lem} \label{lem:  Lemma of zeros of Lazy RW (new)}
Let $S(\cdot)$ be a lazy random walk on $\Z$, starting at
$S(0)=0$, with holding probability $1-\alpha$.  That is,
$$ S(n) = \sum_{i=1}^n x_i , $$
where $x_i$ are i.i.d., such that $\Pr \br{ x_i = 0 } = 1-\alpha$,
and
$$ \Pr \br{ x_i = 1 } = \Pr \br{ x_i = -1 } =
\frac{\alpha}{2} .$$ %

Let $L(n)$ be the number of times the walk visits $0$ up to time
$n$. That is,
$$ L(n) = \abs{ \Set{ 1 \leq i \leq n }{ S(i) = 0 } } . $$
Then, for any $\eps > 0$ there exists $C = C(\eps,\alpha) > 0$
such that for all $n \geq 1$,
\begin{eqnarray*}
    \Pr \br{ L \sr{ \lceil C n^2 \rceil } < n } \leq \eps .
\end{eqnarray*}
\end{lem}

\begin{proof}
Let $m(n)$ be the number of times the walk moves in the first $n$
steps.  Then, $m(n) = \sum_{i=1}^n r_i$, where $r_i$ are i.i.d.
Bernoulli random variables of mean $\alpha$.  By the Chernoff
bound (see e.g. Appendix A in \cite{ProbMethod}),
$$ \Pr \br{ m(n) \leq \frac{\alpha}{2} n } \leq \Pr \br{ \abs{
m(n) - \alpha n} \geq \frac{\alpha}{2} n } < 2 \exp \sr{ -
\frac{\alpha^2}{2} n } . $$

Conditioned on $m(n)$, the walk is a $m(n)$-step simple random
walk.  Thus, for $2k \leq m$, by Corollary \ref{cor: Stirling on
L(n)},
$$ \Pr \Br{ L(n) < k }{ m(n) = m } \leq \frac{k}{\sqrt{m - 2k} } .
$$ %
Let $C > \frac{4}{\alpha}$ and set $j = \lceil C n^2 \rceil$.  If
$i \geq \frac{\alpha j}{2}$ then $2n \leq i$.  Thus,
\begin{eqnarray*}
    \Pr \br{ L(j) < n } & \leq & \Pr \br{ m(j) \leq
    \frac{\alpha}{2} j } + \sum_{i \geq (\alpha j)/2 } \Pr \br{
    L(j) < n \ , \ m(j) = i } \\
    & \leq & \exp \sr{ - \frac{\alpha^2}{2} C n^2 } +
    \frac{\sqrt{2} n}{\sqrt{ \alpha C n^2 - 4 n} } .
\end{eqnarray*}
For large enough $C$ this is less than $\eps$.
\end{proof}

\begin{lem} \label{lem:  Lazy RW does not pass x in x^2 steps}
Let $S(\cdot)$ be a lazy random walk on $\Z$, starting at
$S(0)=0$, with holding probability $1-\alpha$.  That is,
$$ S(n) = \sum_{i=1}^n x_i , $$
where $x_i$ are i.i.d., such that $\Pr \br{ x_i = 0 } = 1-\alpha$,
and
$$ \Pr \br{ x_i = 1 } = \Pr \br{ x_i = -1 } =
\frac{\alpha}{2} .$$ %

Let $m \geq 1$.  Then, for all $\beta > 0$,
$$ \Pr \br{ \max_{1 \leq i \leq m} \abs{ S(i)} < \sqrt{ \beta \alpha m } }
\geq 1 - \frac{1}{\beta} . $$
\end{lem}

\begin{proof}
The assertion is equivalent to
$$ \Pr \br{ \max_{1 \leq i \leq m} \abs{ S(i) } \geq \sqrt{\beta
\alpha m } } \leq \frac{1}{\beta} . $$ %
But this follows immediately from the Kolmogorov inequality, since
$S(i)$ is the sum of i.i.d. random variables, and $\Var \br{ S(m)
} = \alpha m$.
\end{proof}

\section{Random walks on finite graphs}

In this section we recall some properties of a simple random walk
on a finite graph.

Given a finite $d$-regular graph $G$ we define two matrices, whose
columns and rows are indexed by the vertices of the graph.  The
\emph{adjacency matrix} of $G$ is the matrix $A(u,v) = \1{u \sim
v}$ for all $u,v \in G$.  The \emph{transition matrix} of $G$ is
the matrix $P = \frac{1}{d} A$.  It is well known that the
eigenvalues of $P$ are all real.  Further, if $\lambda_1 \geq
\lambda_2 \geq \cdots \geq \lambda_{\abs{G}}$ are the eigenvalues
of $P$, then $\lambda_1 = 1$, and if $G$ is not bi-partite
$\abs{\lambda_i} < 1$ for $1 < i \leq \abs{G}$.  We denote by
$\lambda = \max_{i > 1} \abs{\lambda_i}$.  $\lambda$ is called the
\emph{second eigenvalue} of $G$, and $1-\lambda$ is called the
\emph{spectral gap}.

The following lemma is standard in the theory of random walks on
graphs, and in fact stronger statements can be proved.  We omit
the proof (see \cite{ProbMethod}).

\begin{lem} \label{lem:  mixing time of RW}
Let $G$ be a non-bi-partite $d$-regular graph.  Let $\lambda$ be
the second eigenvalue of $G$.  Let $\mu(i), i \in G$ be \emph{any}
distribution on the vertices of $G$.  Let $x_0,x_1, \ldots, x_t$
be a random walk on $G$, such that $x_0$ is distributed like
$\mu$.  Then, for any $j \in G$,
$$ \abs{ \Pr \br{ x_t = j } - \frac{1}{n} } \leq \lambda^t . $$
\end{lem}

We now prove that the spectral gap of a graph, measures how close
the random walk on the graph is to independent sampling of the
vertices. This is a slight generalization of results from Chapter
9 of \cite{ProbMethod}, and the proof is similar.

In what follows $G$ is a $d$-regular graph of size $n$.  $A$ is
its adjacency matrix. $\lambda$ is the second eigenvalue of the
\emph{transition matrix}.  Thus, $d$ is the largest eigenvalue of
$A$, and all other eigenvalues are at most $d\lambda$.

Let $C \subseteq V(G)$ of size $\abs{C} = cn$.  Define the matrix
$$ Q_C(i,j) = \left\{ \begin{array}{lr}
                            A(i,j) & \textrm{ if } j \in C, \\
                            0 & \textrm{ otherwise. }
                \end{array} \right. $$

For two vectors we use the usual inner product $\IP{ x,y} =
\sum_{i} x(i) \ov{y(i)}$, and norm $\norm{x}^2 = \IP{x,x}$.

\begin{clm} \label{clm:  norm of projection matrix}
$$ \norm{Q_C} \leq \sqrt{ cd^2 + (1-c) d^2 \lambda^2 } $$
\end{clm}

\begin{proof}
Let $x$ be any vector, and let $\tilde{x}$ be the vector defined
by
$$ \tilde{x}(i) = \left\{ \begin{array}{lr}
                            x(i) & \textrm{ if } i \in C, \\
                            0 & \textrm{ otherwise. }
                \end{array} \right. $$
Then $Q_C x = Q_C \tilde{x} = A \tilde{x}$.  Also, note that
$$ \norm{x}^2 = \sum_i x(i)^2 \geq \sum_{i \in C} x(i)^2 =
\norm{\tilde{x}}^2 . $$ %
Thus,
$$ \norm{Q_C}^2 = \max_{x \neq 0} \frac{ \IP{ Q_Cx, Q_Cx}
}{ \IP{x,x} } \leq \max_{\tilde{x} \neq 0} \frac{ \IP{ A
\tilde{x}, A \tilde{x} } }{ \IP{\tilde{x}, \tilde{x} } } . $$ %
So it is enough to prove that for all $x$ such that $\norm{x} = 1$
and such that $x(i) = 0$ for all $i \not\in C$, that $\IP{ Ax, Ax
} \leq cd^2 + (1-c) d^2 \lambda^2$.  Let $x$ be a vector such that
$x(i) = 0$ for all $i \not\in C$, and assume that $\norm{x} = 1$.
Let $\gamma_1 \geq \gamma_2 \geq \ldots \geq \gamma_n$ be the
eigenvalues of $A$, and let $v_1,\ldots,v_n$ be the orthonormal
basis of eigenvectors of $A$, corresponding to these eigenvalues.
We have that $v_1 = n^{-1/2} e$ where $e$ is the all-ones vector.
Decompose $x$,
$$ x = \sum_{s=1}^n \alpha_s v_s . $$
So, by the Cauchy-Schwartz inequality,
$$ \alpha_1 = \IP{x,v_1} = \sum_{i \in C} x(i) \frac{1}{\sqrt{n}} \leq
\sqrt{ \sum_{i \in C} x(i)^2 } \cdot \sqrt{ \sum_{i \in C}
\frac{1}{n} } = \sqrt{c} . $$ %
Note that $\sum_s \alpha_s^2 = \norm{x} = 1$. Thus,
$$ \IP{Ax,Ax} = \sum_{s=1}^n \gamma_s^2 \alpha_s^2 \leq d^2
\alpha_1^2 + (1-\alpha_1^2) (d \lambda)^2 \leq c d^2 + (1-c) (d
\lambda)^2 . $$
\end{proof}

\begin{clm} \label{clm:  Lambda bound}
Let $C_1,C_2,\ldots,C_t$ be subsets of $V(G)$ such that $\abs{C_s}
= c_sn$ for all $s$. Let $\Lambda$ be the number of paths
$x_0,x_1,\ldots,x_t$ in $G$ such that $x_s \in C_s$ for all $s
\geq 1$. Then,
$$ \Lambda \leq n \prod_{s=1}^t \sqrt{ c_s
d^2 + (1-c_s) d^2 \lambda^2 } . $$
\end{clm}

\begin{proof}
For $1 \leq s \leq t$, let $Q_s = Q_{C_s}$.  Let $Q = Q_1 Q_2
\cdots Q_t$.  We claim that
\begin{equation} \label{eqn:  Q(i,j)}
Q(i,j)  \textrm{ is the number of paths } i=x_0,x_1, \ldots,x_t=j
\textrm{ such that } x_s \in C_s \textrm{ for all } s \geq 1.
\end{equation}

This is proven by induction on $t$.  For $t=1$, $Q = Q_1$.  So
$Q(i,j) = 1$ iff $j \in C_1$ and $i \sim j$, and $Q(i,j) = 0$
otherwise. Assume (\ref{eqn:  Q(i,j)}) for $t-1$.  Let $Q' = Q_1
Q_2 \cdots Q_{t-1}$.  Then by the induction hypothesis, $Q'(i,k)$
is the number of paths $i = x_0,x_1,\ldots, x_{t-1} = k$ such that
$x_s \in C_s$ for all $1 \leq s \leq t-1$.  Thus,
$$ Q(i,j) = (Q'Q_t)(i,j) = \sum_k Q'(i,k) Q_t(k,j) $$
is the required quantity.

Thus, if $e$ is the all-ones vector, using (\ref{eqn:  Q(i,j)})
and claim \ref{clm:  norm of projection matrix}, we get that
\begin{eqnarray*}
    \Lambda & = & \sum_{i,j} Q(i,j) = \IP{ Q e, e}  \leq \IP{e,e} \norm{Q} \\
    & \leq & n \prod_{s=1}^t \norm{Q_s} \leq n
    \prod_{s=1}^t \sqrt{ c_s d^2 + (1-c_s) d^2 \lambda^2 }
\end{eqnarray*}
\end{proof}

\begin{lem} \label{lem:  RW stays in specific sets}
Let $x_0,x_1,\ldots,x_t$ be a random walk on $G$ starting at a
uniformly chosen vertex.  Let $C_1,C_2,\ldots,C_t$ be subsets of
$V(G)$ such that $\abs{C_s} = c_sn$ for all $s$.  Let $E$ be the
event that $x_s \in C_s$ for all $s \geq 1$.  Set $c = \sum_s
(1-c_s)$. Then,
$$ \Pr [E] \leq \exp \sr{- \frac{c}{2} (1-\lambda) } . $$
\end{lem}

\begin{proof}
The total number of possible paths is $nd^t$.  Thus, by claim
\ref{clm:  Lambda bound},
\begin{eqnarray*}
    \Pr[E] & = & \frac{\Lambda}{nd^t} \leq \prod_{s=1}^t \sqrt{ c_s +
    (1-c_s) \lambda^2 } \\
    & = & \prod_{s=1}^t \sqrt{ 1 - (1-c_s) + (1-c_s)
    \lambda^2 } \leq \prod_{s=1}^t \exp \sr{ - \frac{(1-c_s)}{2}
    (1-\lambda^2) } \\
    & = & \exp \sr{ - \sum_{s=1}^t \frac{(1-c_s)}{2} (1-\lambda^2) }
    < \exp \sr{- \frac{c}{2} (1-\lambda) } .
\end{eqnarray*}
\end{proof}




\end{document}